\documentclass[11pt]{article}
\usepackage[left=.75in,right=.75in,top=.75in,bottom=.75in]{geometry}
\usepackage{placeins,float,fewerfloatpages,microtype}
\usepackage{titlesec}
\titleformat{\subsection}{\normalfont\normalsize\bfseries}{\thesubsection}{.7em}{}
\titlespacing*{\subsection}{0pt}{9pt plus 2pt minus 1pt}{4pt plus 1pt}
\usepackage{caption}
\AtBeginDocument{%
\setlength{\abovedisplayskip}{7pt plus 2pt minus 2pt}%
\setlength{\belowdisplayskip}{7pt plus 2pt minus 2pt}%
\setlength{\abovedisplayshortskip}{4pt plus 2pt}%
\setlength{\belowdisplayshortskip}{7pt plus 2pt minus 2pt}}
\usepackage{amsmath,amssymb,amsthm,booktabs,array,graphicx,tikz,hyperref,needspace}
\usetikzlibrary{arrows.meta,calc,positioning,decorations.markings}
\hypersetup{colorlinks=true,linkcolor=linkgreen,citecolor=blue,urlcolor=blue,pdfauthor={Anar Akhmedov},pdftitle={Braided Surfaces and Akhmedov's Fake Quadrics}}
\definecolor{myred}{RGB}{205,55,68}
\definecolor{myblue}{RGB}{45,95,185}
\definecolor{mygreen}{RGB}{30,145,100}
\definecolor{linkgreen}{RGB}{0,85,35}
\definecolor{myturq}{RGB}{0,150,170}
\definecolor{mypink}{RGB}{200,70,150}
\definecolor{mypurple}{RGB}{115,70,175}
\newcommand{\Z}{\mathbb Z}\newcommand{\Q}{\mathbb Q}
\newcommand{\Int}{\operatorname{int}}
\newcommand{\coker}{\operatorname{coker}}

\newtheorem{theorem}{Theorem}[section]\newtheorem{proposition}[theorem]{Proposition}\newtheorem{corollary}[theorem]{Corollary}
\newcommand{\fiberpicture}[1]{%
 \begin{tikzpicture}[x=0.65cm,y=0.67cm,>=Stealth,font=\small]
  \pgfmathsetmacro{\height}{0.5*(#1-1)+0.48}
  \foreach \j in {0,...,\numexpr#1-1\relax}{
   \pgfmathsetmacro{\y}{0.5*(#1-1)-\j}
   \ifcase\j\def\col{myred}\or\def\col{myblue}\or\def\col{mygreen}\or\def\col{mypurple}\fi
   \filldraw[fill=\col!20,draw=\col,thick] (-1.4,\y-0.16) rectangle (1.4,\y+0.16);
   \draw[\col,very thick,->] (-0.65,\y)--(0.65,\y);
   \node[right,text=\col] at (2.12,\y) {$c_{\j}$};
  }
  \filldraw[fill=myturq!15,draw=myturq,thick] (-1.4,0) ellipse [x radius=0.70,y radius=\height];
  \filldraw[fill=myturq!15,draw=myturq,thick] (1.4,0) ellipse [x radius=0.70,y radius=\height];
  \node at (-1.4,0) {$N$};\node at (1.4,0) {$S$};
  \foreach \j in {1,...,\numexpr#1-1\relax}{
   \pgfmathsetmacro{\gap}{0.5*(#1-1)-\j+0.5}
   \node at (0,\gap-0.02) {$x_{\j}$};
  }
  \node[below] at (0,-\height-0.1) {$x_{#1}$: outer boundary};
  \node[above,font=\bfseries] at (0,\height+0.3) {$P_{#1}$: two disks, $#1$ bands};
 \end{tikzpicture}}
\newcommand{\basepolygon}[1]{%
 \begin{tikzpicture}[scale=0.66,font=\scriptsize,>=Stealth]
 \pgfmathtruncatemacro{\n}{4*#1}
 \foreach \k in {0,...,\numexpr4*#1-1\relax}{
  \pgfmathsetmacro{\ang}{90-360*\k/\n}
  \coordinate (v\k) at (\ang:1.8);
 }
 \draw[fill=myturq!6] (v0)\foreach \k in {1,...,\numexpr4*#1-1\relax}{--(v\k)}--cycle;
 \foreach \k in {0,...,\numexpr4*#1-1\relax}{
  \pgfmathtruncatemacro{\next}{mod(\k+1,\n)}
  \pgfmathtruncatemacro{\i}{floor(\k/4)+1}
  \pgfmathtruncatemacro{\j}{mod(\k,4)}
  \ifcase\j\def\lab{a_{\i}}\or\def\lab{b_{\i}}\or\def\lab{a_{\i}^{-1}}\or\def\lab{b_{\i}^{-1}}\fi
  \ifnum\k=0\def\col{myred}\else\ifnum\k=2\def\col{myred}\else\def\col{myblue}\fi\fi
  \draw[\col,thick] (v\k)--(v\next);
  \pgfmathsetmacro{\mid}{90-360*(\k+0.5)/\n}
  \node[text=\col] at (\mid:2.14) {$\lab$};
 }
 \node[align=center] at (0,0) {$D^2$\\base $\Sigma_{#1}$};
 \end{tikzpicture}}
\newcommand{\surfacepairs}[1]{%
 \begin{tikzpicture}[x=1.23cm,y=0.72cm,font=\small,>=Stealth]
  \filldraw[fill=myturq!7,draw=myturq] (-0.75,-0.72) rectangle (#1-0.25,0.65);
  \foreach \j in {0,...,\numexpr#1-1\relax}{
   \draw[fill=white,draw=myturq,thick] (\j,0) circle (0.36);
   \draw[mypurple,thick] (\j-0.3,0.2) .. controls (\j-0.8,1.7) and (\j+0.8,1.7) .. (\j+0.3,0.2);
   \draw[myred,thick] (\j-0.20,0.65) .. controls (\j-0.08,1.28) and (\j+0.08,1.28) .. (\j+0.20,0.65);
   \draw[myblue,thick] (\j-0.28,0.85) .. controls (\j,0.58) .. (\j+0.28,0.85);
   \node at (\j,0) {$T_{\j}$};
   \node[above] at (\j,1.38) {$(e_{\j},f_{\j})$};
  }
  \node[below] at (#1/2-0.5,-0.78) {planar core; outer circle capped by a disk};
 \end{tikzpicture}}

\title{Braided Multisections and Akhmedov's Fake Quadrics}
\author{Anar Akhmedov}
\date{October 28, 2026}
\begin{document}\maketitle
\begin{abstract}
We give handlebody decompositions of the complements of the cyclic
multisections in $\Sigma_g\times S^2$ constructed in~\cite{AB}.
The 2-handle framings are induced by the product coordinates.
We compute the integral homology of the complements from their
handle chain complexes, including the torsion subgroup of $H_2$.
For the degree pairs $(3,3)$ and $(2,4)$, we specify surface
diffeomorphisms inducing the maps on first homology prescribed
in~\cite{AY} and recover the integral first homology groups
$(\Z/3)^2\oplus\Z/(3m)$ and $(\Z/2)^2\oplus\Z/(4m)$,
respectively, for $m\geq1$. For the symmetric $m=1$ sum, we also give
a simultaneous surface diagram of the relative product attachments
and identify their product framings with the surface framings.
The resulting symplectic sums are
aspherical and symplectically minimal, and have the rational cohomology
ring of $S^2\times S^2$.
\end{abstract}

\section{Introduction}
The manifolds constructed in~\cite{AB} have the rational cohomology
ring of $S^2\times S^2$. We call them Akhmedov's fake quadrics;
no complex structure is assumed. They are symplectic sums of ruled
surfaces along cyclic multisections. We give handlebody decompositions
of their complements, compute the integral handle chain complexes,
identify the $H_2$-torsion geometrically, and specify marked surface
gluing maps. For the symmetric $(3,3)$ sum with twist parameter
$m=1$, we also give a simultaneous surface diagram of the relative
product attachments and determine their surface framings.

The earlier construction in~\cite{A06} used knot surgery and twisted
fiber sums to obtain the integral cohomology ring of $S^2\times S^2$.
The cyclic construction in~\cite{AB} instead uses multisection
complements. It determines their groups and boundary incompressibility,
classifies the degree pairs with Euler characteristic four and signature
zero, and gives a rational-homology gluing criterion. The resulting
infinite families are distinguished by the orders of their finite first
homology groups. For related non-isotopic symplectic surfaces, see
Smith~\cite{Smith} and Hays and Park~\cite{HP}.

In joint work with Sai-Kee Yeung~\cite{AY}, we computed integral
Mayer-Vietoris presentations and first-homology formulas for the three
product families. Here we realize the matrices of~\cite[Theorem~4.5]{AY}
by marked surface diffeomorphisms. Their based mapping classes are part
of the gluing data: the matrices alone do not determine the smooth sums.

Section~\ref{sec:complements} gives the complement handles and homology;
Sections~\ref{sec:markings} and~\ref{sec:families} give the markings and
gluings. The simultaneous surface-framed relative attachments are in
Theorem~\ref{thm:surfacecarrier}. In the follow-up paper with Selman
Akbulut and Eylem Zeliha Yildiz~\cite{AAyII}, we will study full
handlebody diagrams of the closed sums, including the complete attaching
links and their integer framing coefficients.

Removing a regular sphere fiber as well gives the Stein fillings
of~\cite{Stein}. For fixed $g\geq2$ and $p\geq4$, these include
infinitely many pairwise nonhomeomorphic aspherical fillings of one
contact three-manifold, with the same integral homology and zero
intersection pairings. Their groups are distinguished by the first
homology of finite covers. For related boundary gluings and handlebody
constructions, see~\cite[Chapter~3 and Section~4.3]{AkBook}
and~\cite{AkCCM,AkAP3,AkAP2}.

Throughout, $\Sigma_g$ is a closed oriented surface of genus
$g\geq2$, and $p\geq2$ is the multisection degree. Identify the
oriented equator of $S^2$ with $S^1=\{z\in\mathbb C:|z|=1\}$.
Write $\mathcal S=\Sigma_g\times\{\mathrm{pt}\}$ and choose a smooth
map $u:\Sigma_g\to S^1$ for which
$u_*:\pi_1(\Sigma_g)\to\mathbb Z$ is surjective. Set
\[
 F_{g,p}=\{(x,z)\in\Sigma_g\times S^1:z^p=u(x)\}\subset\Sigma_g\times S^2,
 \qquad M_{g,p}=(\Sigma_g\times S^2)\setminus\Int\nu F_{g,p}.
\]
Here $\nu F_{g,p}$ is a closed tubular neighborhood. The first projection
$F_{g,p}\to\Sigma_g$ is a connected, unbranched $p$-fold cover;
connectedness follows from the surjectivity of $u_*$ modulo $p$.
The product symplectic form restricts to the pullback of the base
area form, since the sphere coordinate lies in the equator.
The tubular disks are chosen fiberwise and equivariantly. The multisection
has class $p[\mathcal S]$, square zero, and genus $G=1+p(g-1)$. The complement
is a bundle with fiber $P_p$, a sphere with $p$ disks removed.

\begin{figure}[tbp]
\centering
\includegraphics[width=.98\textwidth,height=.58\textheight,keepaspectratio]{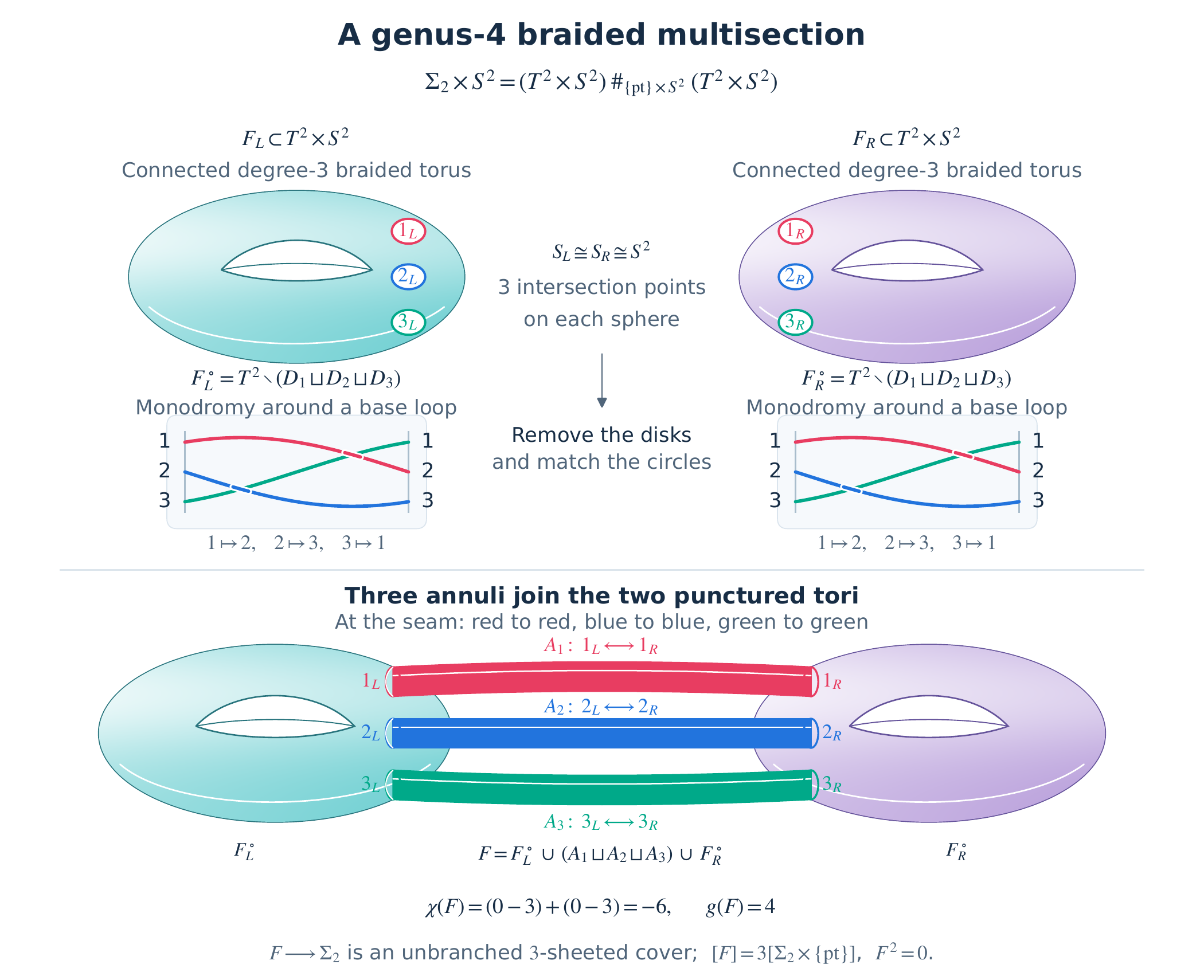}
\caption{The genus-four braided multisection in $\Sigma_2\times S^2$.
Each braided torus meets the summing sphere $\{\mathrm{pt}\}\times S^2$
in three points. The punctured tori are joined by three annuli, giving
$\chi(F)=-3-3=-6$ and $g(F)=4$. At the seam, $j_L$ is joined to $j_R$, with matching colors.
The large silhouettes show the surface pieces and their annular joins;
the insets record $1\mapsto2\mapsto3\mapsto1$ around a base loop.
Strand colors follow moving points, while numbers label fixed fiber positions.}
\label{fig:genus4-braided-sphere-sum}
\end{figure}

Figure~\ref{fig:genus4-braided-sphere-sum} gives the $g=2$, $p=3$
model as a sum of two connected degree-three braided tori in
$T^2\times S^2$. Here $u$ represents the primitive class $(1,0,1,0)$.
An orientation-preserving change of base marking gives the fixed
convention used in Section~\ref{sec:complements}.
Take the two defining maps constant near the summing points, with the
same value. Choose product tubular neighborhoods $S^2\times D^2$ of
the summing spheres, with identical three-point configurations and deleted disks.
On their boundaries the gluing is
\[
 \Psi(v,e^{i\theta})=(v,e^{-i\theta}),\qquad v\in S^2.
\]
Thus $j_L$ joins $j_R$ for $j=1,2,3$ through a product annulus.
The cyclic braid occurs around a base loop in each torus, with no
additional braid at the sphere-sum seam.

The product cases of~\cite{AB} are listed below; $G$ is the common
gluing genus.

\begin{center}
\begin{tabular}{cccc}
\toprule Degrees & Pieces & $G$ & $(\chi(M_L),\chi(M_R))$\\\midrule
$(2,3)$ & $M_{4,2},\ M_{3,3}$ & $7$ & $(0,4)$\\
$(2,4)$ & $M_{3,2},\ M_{2,4}$ & $5$ & $(0,4)$\\
$(3,3)$ & $M_{2,3},\ M_{2,3}$ & $4$ & $(2,2)$\\\bottomrule
\end{tabular}
\end{center}
We study the $(3,3)$ and $(2,4)$ families in product ruled surfaces.

Theorem~\ref{thm:piecehomology} gives
\[
 H_2(M_{g,p};\Z)\cong\Z^{1+2(g-1)(p-1)}\oplus\Z/p,
\]
and Corollary~\ref{cor:poletorsion} identifies the torsion generator.
The closed first-homology formulas in the following theorem recover
those of~\cite[Theorem~4.5]{AY}.

\begin{theorem}\label{thm:main}
The complements $M_{2,3}$, $M_{3,2}$, and $M_{2,4}$ admit
the handlebody decompositions described in Section~\ref{sec:complements},
with respective handle counts $(h_0,h_1,h_2,h_3)$ given by
\[
 (2,11,14,3),\qquad(2,14,14,2),\qquad(2,12,18,4).
\]
For each integer $m\geq1$, let $Z_m^{3,3}$ and $Z_m^{2,4}$ be the symplectic sums
formed using the marked surface diffeomorphisms $f_m$ of
Sections~\ref{sec:markings} and~\ref{sec:families}. In the chosen boundary
trivializations, the gluing is
\[
 \phi_m(w,e^{i\theta})=(f_m(w),e^{-i\theta}),\qquad w\in\Sigma_G,
\]
with no surface-dependent normal-circle twist. Then
\[
\begin{aligned}
H_1(Z_m^{3,3};\Z)&\cong(\Z/3)^2\oplus\Z/(3m),\\
H_1(Z_m^{2,4};\Z)&\cong(\Z/2)^2\oplus\Z/(4m).
\end{aligned}
\]
Both symplectic sums have Euler characteristic four, signature zero, and the
rational cohomology ring of $S^2\times S^2$. They are aspherical and
symplectically minimal.
\end{theorem}
\Needspace{7\baselineskip}
\medskip
\noindent\textbf{Acknowledgments.}
The author is very grateful to B. Doug Park for several valuable
discussions in 2007--2010 concerning one of the fake quadrics studied here,
and for his insights and encouragement. The author also thanks Selman
Akbulut, Eylem Zeliha Yildiz, Sai-Kee Yeung, \c{C}a\u{g}r\i{} Karakurt, and
S\"umeyra Sakall\i{} for helpful discussions. An LLM-based tool assisted
with some algebraic computations involving the homology presentations
and gluing matrices, preparation of the computational scripts and
interval-arithmetic verification described in Appendix~\ref{app:crossings},
and figure preparation. All mathematical ideas and results are due to the
author, who takes full responsibility for the content.

\medskip
\noindent\textbf{Computational materials.}
The scripts, input data, surface-routing and cap-tree data, and interval
certificates used in this paper are available in the
\href{https://github.com/akhmedov-coder/fake-quadrics-computations}{accompanying GitHub repository}.
The file \path{COMPUTATIONAL_README.txt} gives reproduction commands,
dependencies, software versions, and the scope of each check.
The homology computations are also proved in the text; the accompanying
files provide the finite data used in the computer-assisted diagram
verifications.

\section{Handlebody decompositions of the complements}\label{sec:complements}
Products of paths are traversed from left to right.
Write $R_\theta$ for rotation about the polar axis through angle
$\theta$, and put $r_p=R_{2\pi/p}$. Choose based generators
$a_i,b_i$ of $\pi_1(\Sigma_g)$ with relation
$\prod_{i=1}^g[a_i,b_i]=1$, where $[a,b]=aba^{-1}b^{-1}$.
We use the same letters for the corresponding symplectic homology
basis, with $a_i\cdot b_j=\delta_{ij}$.
Use a marking with $u_*(a_1)=-1$, $u_*(b_1)=0$, and
$u_*(a_i)=u_*(b_i)=0$ for $2\leq i\leq g$, and write
$\chi_p=-u_*\bmod p:\pi_1(\Sigma_g)\to\mathbb Z/p$. Such a marking exists since $u_*$ is surjective.
Our face-pairing convention uses $r_p$ for $a_1$ and the identity
for the other generators. Positive geometric transport is the inverse
of this quotient action. Conjugation by the displayed $a_1$ acts
by $r_p$ on the fiber group. Positive passage through the corresponding
1-handle has the inverse angular shift. For a fiber loop $x$ based at the fixed north pole, the
transport square has boundary
$x a_1(r_p^{-1}x)^{-1}a_1^{-1}$, giving
$a_1xa_1^{-1}=r_p(x)$. This accounts for the inverse in the passage
convention. Let
$\rho:\widetilde\Sigma_g\to\Sigma_g$ be the universal covering;
$\gamma\in\pi_1(\Sigma_g)$ acts on it by deck transformations.
With this convention, the bundle is
\[
 M_{g,p}\cong(\widetilde\Sigma_g\times P_p)/\pi_1(\Sigma_g),
 \quad \gamma(x,v)=(\gamma x,r_p^{\chi_p(\gamma)}v),
\]
Equivalently, glue the faces of a base polygon times $P_p$ by these
transition maps.

For the quotient identification, choose a real lift $U$ of the
argument of $u$ on $\widetilde\Sigma_g$.  It satisfies
$U(\gamma x)=U(x)+2\pi u_*(\gamma)$.  Rotation of the fiber by
$U(x)/p$ carries the reference deleted disks to the disks centered
at the roots of $z^p=u(x)$.  The map
$(x,v)\mapsto(\rho(x),R_{U(x)/p}v)$ is invariant under
$(x,v)\mapsto(\gamma x,r_p^{-u_*(\gamma)}v)$, hence descends
to the stated diffeomorphism of complement bundles, carrying their
pole sections and boundary collars with it.

The 1-handles are represented by paired attaching 3-balls. Product
coordinates and their prescribed identifications determine the 2-handle
framings. We do not convert these drawings to dotted-circle notation
with integer coefficients relative to Seifert framings.

\subsection{The fiber and product handle decompositions}
We use $I=[-1,1]$; subscripts distinguish the coordinate factors.
Retain a small disk about each pole of the sphere, denoted $N$ and $S$,
and join them by $p$ disjoint bands $c_0,\ldots,c_{p-1}$ between the
equatorial deleted disks. These form two surface 0-handles and $p$ surface 1-handles. Rotation fixes the two disks setwise and sends $c_j$ to
$c_{j+1}$ cyclically. The rotation preserves the orientation of each band.

\begin{figure}[!htbp]
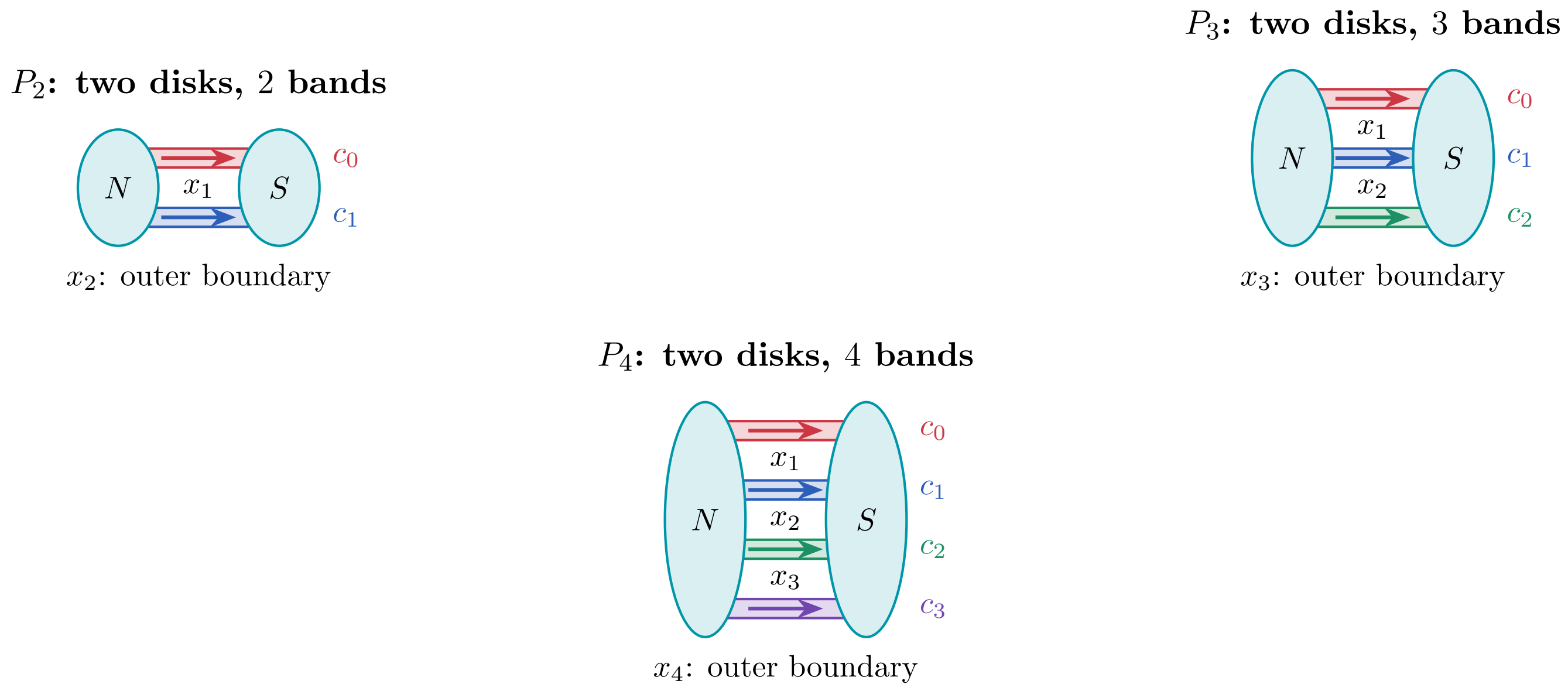

\centering
\fiberpicture{2}\hfill\fiberpicture{3}\par\vspace{8pt}
\fiberpicture{4}
\caption{Planar disk-and-band drawings of the annulus, the
pair of pants, and the four-holed sphere. The turquoise disks are
0-handles; the colored strips are 1-handles. The arrow on each strip is
its core direction from $N$ to $S$. The outer boundary is one of the
$p$ boundary components, not an additional component. Cyclic rotation is
a diffeomorphism of the sphere model; its planar drawing can move an
inner boundary to the outer one.}
\end{figure}

Use $c_j$ also for the oriented core path of the $j$th band,
from its center in $N$ to its center in $S$. Read fiber indices modulo $p$. Define
\[
 x_j=c_{j-1}c_j^{-1}\quad(1\leq j\leq p),\qquad
 x_1\cdots x_p=1,\qquad r_p(x_j)=x_{j+1}.
\]
These are peripheral loops with the indicated base paths at $N$.
After choosing $c_0$ as a maximal tree edge, put
$y_j=c_jc_0^{-1}$, $1\leq j<p$. The free fiber group has the basis
$y_1,\ldots,y_{p-1}$, and
\[
 r_p(y_j)=y_{j+1}y_1^{-1}\ (j<p-1),\qquad
 r_p(y_{p-1})=y_1^{-1}.
\]
The two-disk decomposition retains the cyclic symmetry of the bands,
although $r_p$ does not permute the based free generators.

\par\medskip

For a base band $k\in\{a_i,b_i\}$ and pole $P\in\{N,S\}$,
$k^P$ denotes the 1-handle given by its product with the pole disk.
The fiber band $c_j$ over the base 0-handle gives the connecting
1-handle $C_j$; its core interval is $I_j$.
Write $k^-,k^+$ for the entry and exit feet of the base band. Give $\Sigma_g$ one surface 0-handle, $2g$ surface 1-handles, and one
surface 2-handle. Since $r_p$ preserves the fiber disk-and-band
decomposition, products of base and fiber handles give
four-dimensional handles after rounding corners. The boundary of each
product core attaches through products of lower total index, so the
handles may be ordered by their four-dimensional indices.
Over a base 1-handle
the two disk products attach first, followed by the $p$ band products.
Over the closing base 2-handle, the two disk products attach before the
$p$ band products. The $a_1$ face pairing is $r_p$, so positive passage through the
corresponding 1-handles has the inverse shift; the other base handles use the identity.

\begin{center}
\begin{tabular}{lll}
\toprule Base handle $\times$ fiber handle & Index & Number\\\midrule
$D^2\times N$, $D^2\times S$ & $0$ & $2$\\
base 0-handle $\times c_j$ & $1$ & $p$\\
base 1-handle $\times N$ or $S$ & $1$ & $4g$\\
base 1-handle $\times c_j$ & $2$ & $2gp$\\
base 2-handle $\times N$ or $S$ & $2$ & $2$\\
base 2-handle $\times c_j$ & $3$ & $p$\\\bottomrule
\end{tabular}
\end{center}
The unreduced handle counts are
\[
 (h_0,h_1,h_2,h_3)=(2,\,4g+p,\,2+2gp,\,p).
\]

The product of a base band and a fiber band gives a 2-handle with
core square $[-1,1]_s\times[-1,1]_u$ and normal square
$[-1,1]_t\times[-1,1]_v$. Rounding corners gives its core and normal disks. The $s$-faces attach to the fiber bands
over the endpoints of the base band. The $u$-faces attach to the
base-band products with $N$ and $S$, which have already been added.
The face pairing at $a_1$ identifies $c_j$ with $c_{j+1}$;
positive passage through the corresponding 1-handle therefore joins the $j+1$ entry to the $j$ exit;
along any other base generator it lands on $c_j$. These four face
maps, including their collar coordinates, are the attaching map on
$S^1\times D^2$. The product framing $(\partial_t,\partial_v)$ is transported
by those same maps. Along positive $a_1$ passage the pole-disk map
is $z\mapsto e^{-\sigma_P2\pi is/p}z$, $0\leq s\leq1$,
where $\sigma_N=1$ and $\sigma_S=-1$; its vertical derivative is the same
rotation. Other base 1-handles use the identity. Across $C_j$ the
base angle is fixed and the fiber-band normal matches
$+\partial_{\phi_N}$ with $-\partial_{\phi_S}$, where
$\phi_N,\phi_S$ are angular coordinates in the positively oriented pole disks.
These maps specify the product framing. The flat quotient model
selects the closing trivialization over the base 2-handle. Along each
commutator, the forward and backward rotation paths cancel, with zero
total lifted angular rotation. Thus both the endpoint transport and
the winding in this chosen trivialization are fixed.

\begin{figure}[!htbp]
\centering
\includegraphics[width=.80\textwidth,trim=7.50bp 6.15bp 8.56bp 13.00bp,clip]{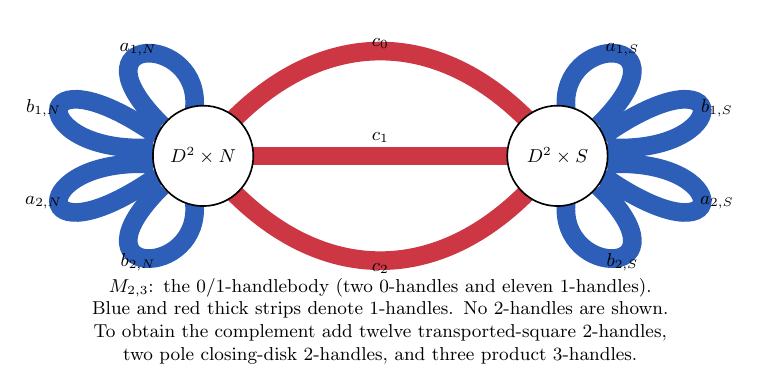}
\caption{The two 0-handles and eleven 1-handles of $M_{2,3}$.
The twelve product 2-handles and two closing 2-handles attach to
this handlebody; their attaching circles are drawn below.}
\label{fig:lower}
\end{figure}

\begin{figure}[!htbp]
\centering
\includegraphics[width=.64\textwidth,trim=57.50bp 5.91bp 57.93bp 6.50bp,clip]{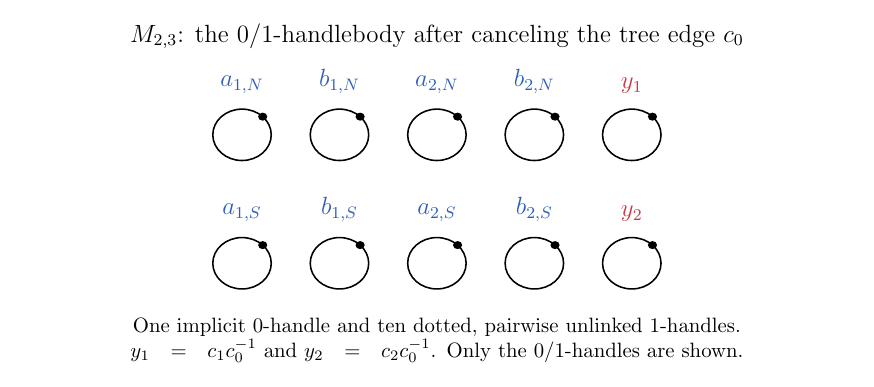}
\caption{The ten 1-handles in dotted-circle notation after canceling
$C_0$ with one 0-handle. Only the 0- and 1-handles are shown.
To express the full attaching link in this notation, the eight strands
through $C_0$ must be carried through the boundary connected-sum
collar. That conversion has not been carried out.}
\label{fig:dottedskeleton}
\end{figure}

\subsection{Attaching maps, coordinate charts, and projections}
Use coordinates $(s,t)$ on a base band and $(u,v)$ on a fiber band.
The normal coordinates of the square core are $(t,v)$. Over each
pole the 0-handle is $D^2_{\rm base}\times D^2_{\rm pole}$.
Its boundary contains the product faces
$\partial D^2_{\rm base}\times D^2_{\rm pole}$ and
$D^2_{\rm base}\times\partial D^2_{\rm pole}$.
Put the base feet on the first face in the order inherited from the
oriented surface polygon and the connecting 1-handle feet on the second face in
the cyclic order of the sphere model. Each square is rounded on
these product faces. Parallel attaching strips retain their product
collars under the face pairings. These collars determine the
attaching curves and their crossings; the relator words alone do not.

In attaching words, $c_j$ denotes passage through $C_j$.
The closing 2-handles are products of the base surface 2-handle
with the pole disks.
With the group convention that $a_1$ conjugates by $r_p$, the oriented attaching circle follows the path
\[
 R_{k,j}=k^N c_j(k^S)^{-1}c_{j+\delta(k)}^{-1},\qquad
 \delta(a_1)=1,\quad\delta(k)=0\ (k\ne a_1).
\]
The label $R_{k,j}$ denotes both the product 2-handle and its
oriented attaching circle; $\delta(k)$ records whether $k=a_1$.
Reversing the base-band convention reverses the shift.
The indices for the two degree pairs are:
\begin{center}
\begin{tabular}{llll}
\toprule Piece & Base generators $k$ & Band indices $j$ & Product 2-handles\\\midrule
$M_{2,3}$ & $a,b,c,d$ & $0,1,2$ & $12$\\
$M_{3,2}$ & $a,b,c,d,e,f$ & $0,1$ & $12$\\
$M_{2,4}$ & $A,B,C,D$ & $0,1,2,3$ & $16$\\\bottomrule
\end{tabular}
\end{center}
The letters in each row abbreviate the ordered base generators
$(a_1,b_1,\ldots,a_g,b_g)$. In particular, for $M_{2,3}$,
$(a,b,c,d)=(a_1,b_1,a_2,b_2)$. Its rotating attaching curves follow
\[
\begin{gathered}
a^Nc_0(a^S)^{-1}c_1^{-1},\qquad a^Nc_1(a^S)^{-1}c_2^{-1},\\
a^Nc_2(a^S)^{-1}c_0^{-1}.
\end{gathered}
\]
For $R_{a,0}$, enter $a^N$ on strand $1$ and leave on strand $0$,
cross $C_0$ to $S$, pass backward through $a^S$ from strand $0$
to strand $1$, and return through $C_1$. The corner charts below
join these passages and transport their normal disks.
The nine nonrotating product 2-handles use matching connecting-handle indices.
For $M_{3,2}$ the two rotating product 2-handles exchange $c_0,c_1$;
for $M_{2,4}$ the four rotating product 2-handles advance cyclically.

\begin{figure}[!htbp]
\centering
\includegraphics[width=.88\textwidth]{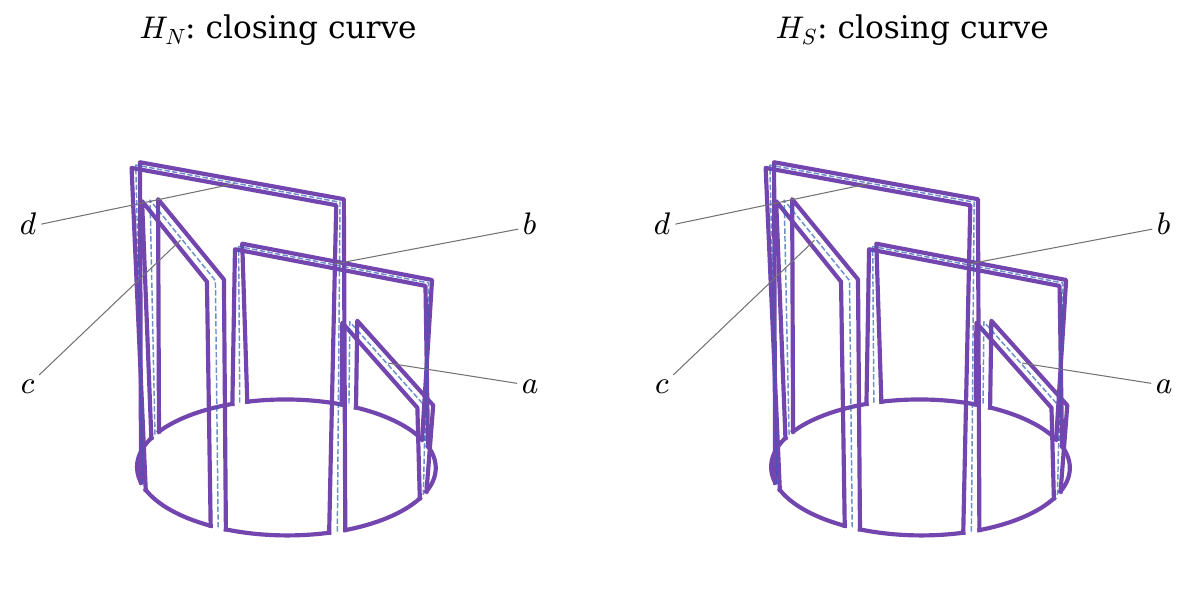}\par
\includegraphics[width=.88\textwidth]{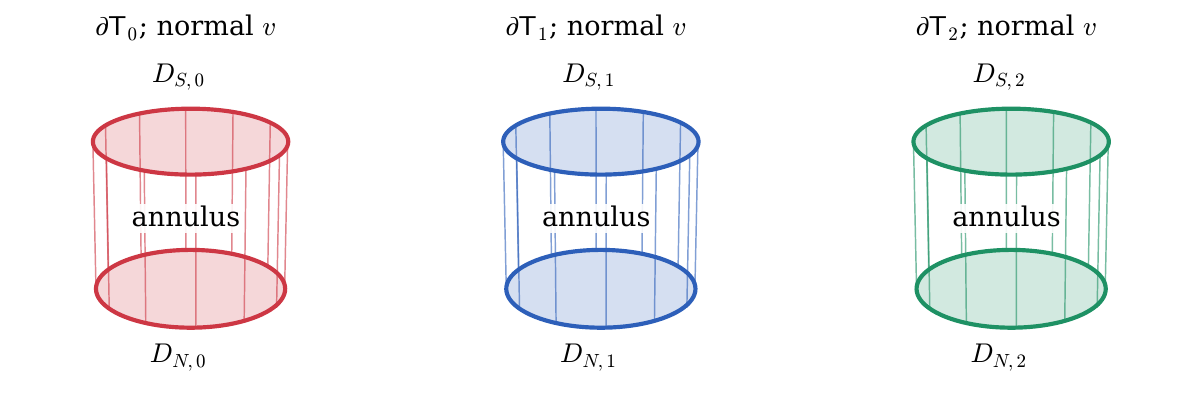}
\caption{The two closing 2-handle attaching curves and the three
3-handle attaching spheres. Each sphere consists of two closing-core
disk faces joined by an annular face in the corresponding fiber-band
collar. The normal interval is the transverse coordinate of that band.
These are local product-face schematics, not a simultaneous planar
projection of the entire post-attachment boundary.}
\label{fig:closingandspheres}
\end{figure}

Denote the two additional 2-handles by $H_N,H_S$ and the $p$
3-handles by $\mathsf T_j$. The 2-handles have base-disk cores at the
pole centers. Their attaching curves are the two base-surface attaching curves
$\prod_i[a_i^N,b_i^N]$ and $\prod_i[a_i^S,b_i^S]$, embedded by the
surface polygon with their pole-disk normal frames. The $p$ 3-handle
attachments are also specified by product faces: for a base closing
handle $D^2_\xi\times D^0$ and a fiber band
$D^1_u\times D^1_v$, their cores are $D^2_\xi\times D^1_u$,
and the attaching sphere is
\[
 (\partial D^2_\xi\times D^1_u)
 \ \cup\ (D^2_\xi\times\partial D^1_u).
\]
The first face meets the previously attached product 2-handles and the
second meets the closing 2-handles $H_N,H_S$. The normal interval is
$D^1_v$, with the transported collar coordinates.
Figure~\ref{fig:closingandspheres} shows these attaching spheres.

\begin{figure}[!htbp]
\centering
\includegraphics[width=.76\textwidth,trim=39.50bp 6.01bp 40.24bp 6.00bp,clip]{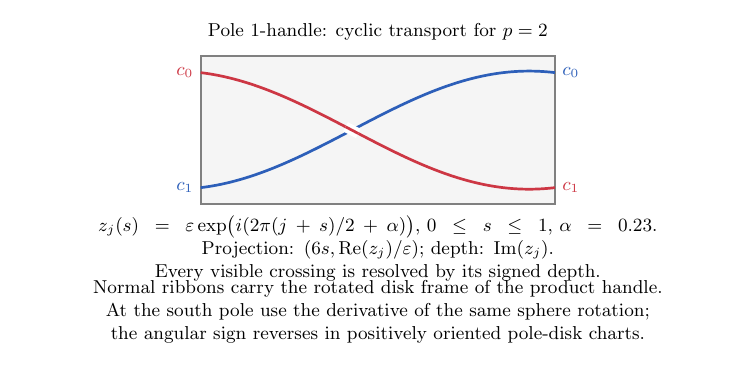}\par
\includegraphics[width=.76\textwidth,trim=39.50bp 6.01bp 40.24bp 6.00bp,clip]{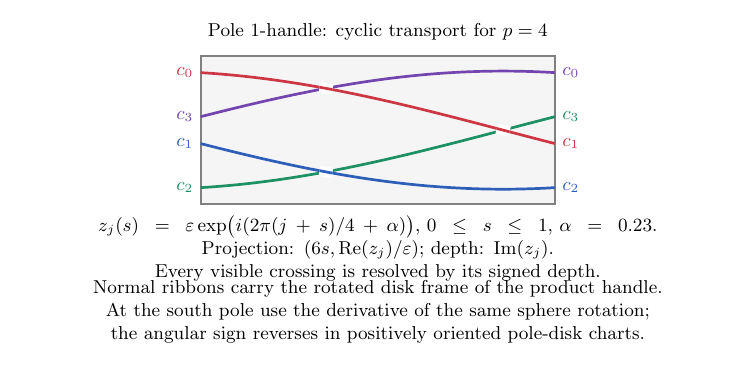}
\caption{The degree-two and degree-four transport maps used in
the mixed case. Here $\varepsilon>0$ is the strand radius in the
pole disk, $j$ labels the fiber strand, and $0\leq s\leq1$.
The northern coordinate is
$z_j(s)=\varepsilon\exp(i(2\pi(j+s)/p+0.23))$.
Here $s$ increases opposite the positive displayed base generator.
These drawings fix the crossings within one pole handle; they are
not the entire attaching link of either complement.}
\label{fig:helix24}
\end{figure}

\begin{center}
\begin{tabular}{crrrrr}
\toprule Piece & $h_0$ & $h_1$ & $h_2$ & $h_3$ & $\chi$\\\midrule
$M_{4,2}$ & 2&18&18&2&0\\
$M_{3,3}$ & 2&15&20&3&4\\
$M_{3,2}$ & 2&14&14&2&0\\
$M_{2,4}$ & 2&12&18&4&4\\
$M_{2,3}$ & 2&11&14&3&2\\\bottomrule
\end{tabular}
\end{center}
The alternating sum is $\sum(-1)^ih_i=(2-2g)(2-p)$.

\par\medskip

We draw the attaching curves in one 0- and 1-handlebody. Order the eight base feet
around each pole by
\[
 a^-,b^+,a^+,b^-,c^-,d^+,c^+,d^-.
\]
Write their distinct base-disk angles as $\theta_{k\pm}$.
Use complex coordinates $(z_{\rm base},z_{\rm pole})$ on each
rounded $S^3\subset\mathbb C^2$, with arguments $\theta,\phi$.
In positively oriented
pole coordinates put
$\phi_{N,j}=2\pi j/3+\alpha$ and
$\phi_{S,j}=-(2\pi j/3+\alpha)$, where $\alpha=0.23$.
The positive constants $\epsilon,\zeta,\eta$ are respectively a
polar cutoff, an angular half-width for the feet, and a normal-disk
scale. Choose disjoint base-foot neighborhoods
$\{|\theta-\theta_{k\pm}|\leq\zeta,\ |z_{\rm pole}|\leq\sin\epsilon\}$
and connecting-handle foot neighborhoods
$\{|\phi-\phi_{P,j}|\leq\zeta,\ |z_{\rm base}|\leq\sin\epsilon\}$.
These are product 3-balls after rounding. The two types are disjoint
because $\epsilon<\pi/4$; within either type take $\zeta$ smaller than
half the foot-angle separation and choose $0<\eta<\zeta$.
For the three pieces, $\zeta=0.05$ gives foot neighborhoods containing
all the attaching strands.
Extend corner endpoints slightly into their rounded collars.

For $\epsilon\leq\tau\leq\pi/2-\epsilon$, use the corner arcs
\[
 \Gamma_{P,k\pm,j}(\tau)=
 \bigl(\cos\tau\,e^{i\theta_{k\pm}},
       \sin\tau\,e^{i\phi_{P,j}}\bigr).
\]
Their local normal-disk charts are obtained by varying the two arguments.
When $(t,v)$ denotes the oriented surface-band coordinates, use
$\theta_{k\pm}+\lambda_{\pm}\eta t$ and
$\phi_{P,j}+\sigma_P\eta v$, where
$\lambda_- =1$, $\lambda_+=-1$, $\sigma_N=1$, and $\sigma_S=-1$.
These signs record the orientation reversal between the two ends of a
surface band. An unsigned local argument chart is equivalent only when its
endpoint identification carries the same reversal. In either notation the
normal frame also rotates along the distinguished base 1-handle.
Choose $\eta$ smaller than half the angular foot separation.
Two different pairs of arguments cannot meet in the interior, so
all twenty-four ribbons at each pole are disjoint. The coordinate atlas uses
$\epsilon=0.18$, $\eta=0.022$; normal disks may be shrunk for the
projections below.

\begin{figure}[!htbp]
\centering
\includegraphics[width=\textwidth]{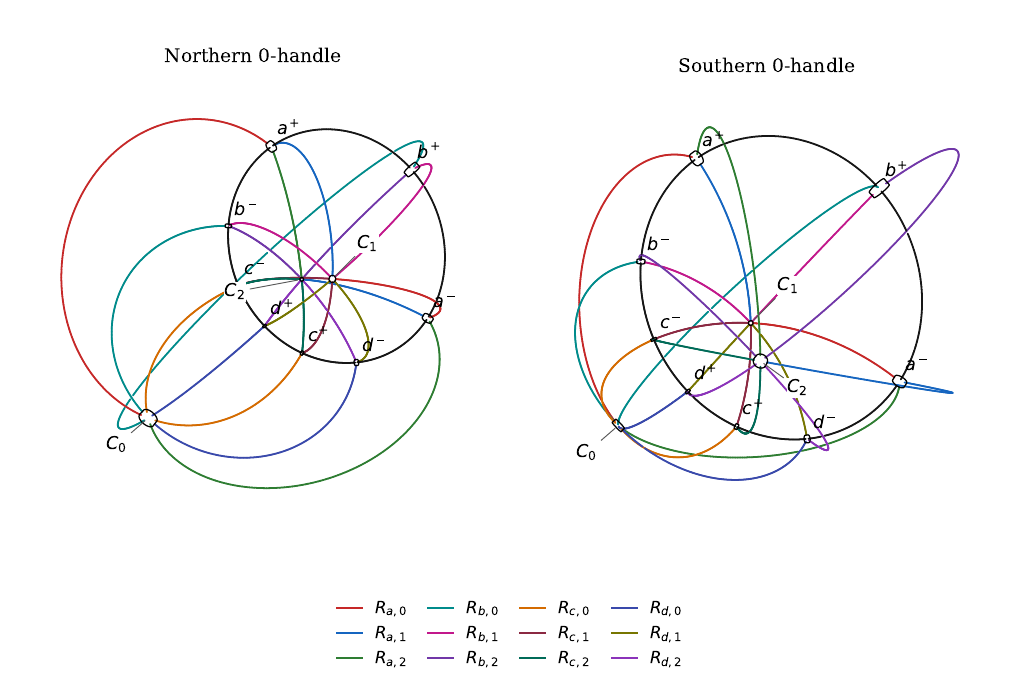}
\caption{$M_{2,3}$. Northern and southern projections of the attaching
curves. Pair $k^-$ with $k^+$ within each pole and $C_j$ across
poles. Colored arcs belong to the product 2-handles and black arcs
to the closing 2-handles. The transparent outlines mark the
1-handle feet. Figures~\ref{fig:usualmatchings}
and~\ref{fig:usualframematchings} give the strand identifications
and product framings.}
\label{fig:corners}\label{fig:usualsouth}
\end{figure}

\begin{figure}[!htbp]
\centering
\includegraphics[width=\textwidth,height=.69\textheight,keepaspectratio,trim=70.00bp 22.40bp 70.00bp 22.50bp,clip]{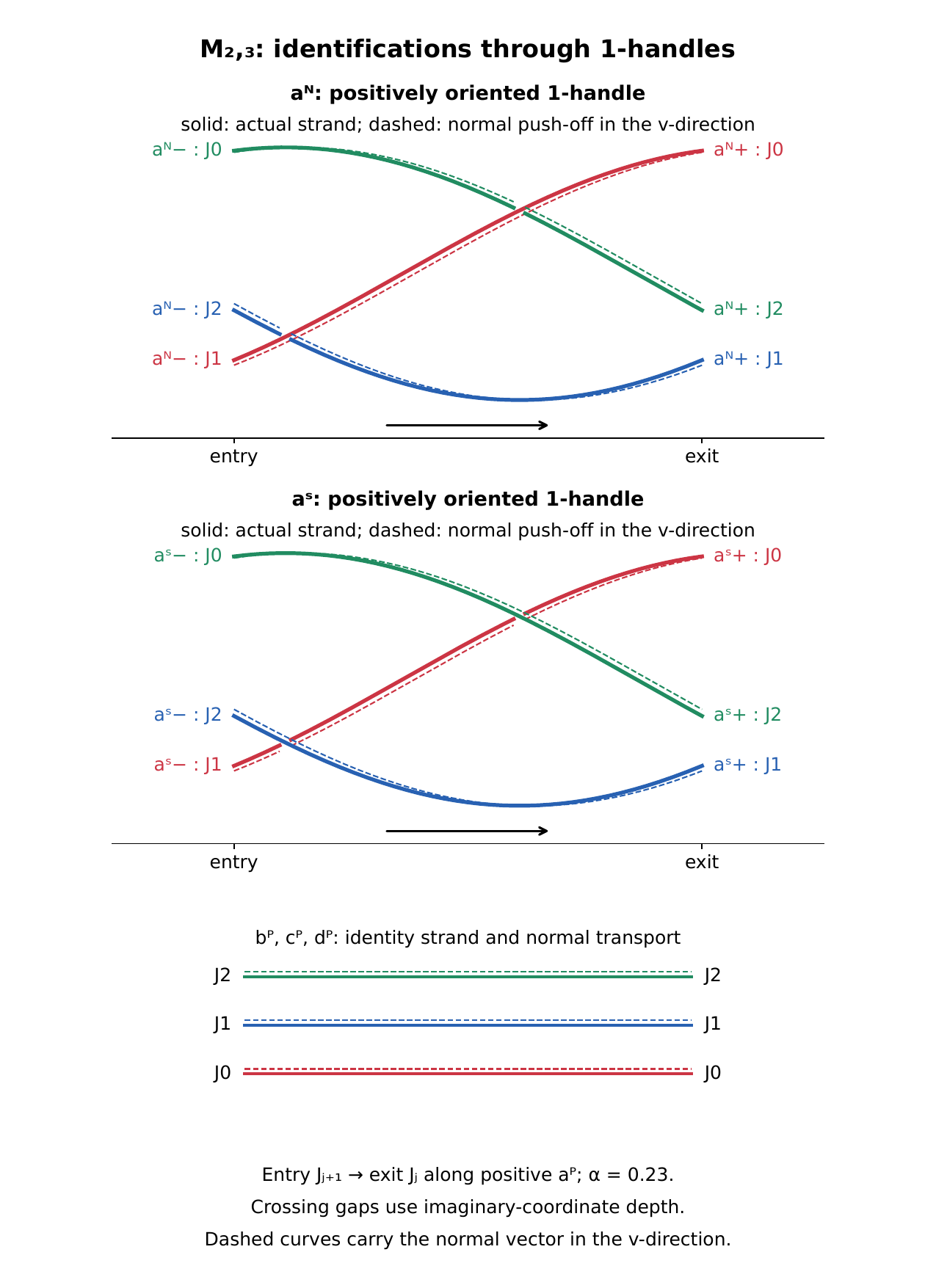}
\caption{Strand identifications through the base 1-handles over $N$ and $S$.
Here $J_j$ labels the position of strand $j$ in a pole-disk cross-section.
The cyclic strands run from $J_{j+1}$ to $J_j$ along positive $a$;
the other base 1-handles use identity transport. Dashed push-offs indicate the fiber-band normal direction. Figure~\ref{fig:usualframematchings} gives the identifications through
the connecting 1-handles and the framing of the closing 2-handle.
Colors label local strands, not globally colored components in
Figure~\ref{fig:corners}.}
\label{fig:usualmatchings}
\end{figure}

\begin{figure}[!htbp]
\centering
\includegraphics[width=\textwidth,height=.69\textheight,keepaspectratio,trim=47.00bp 28.40bp 51.50bp 22.50bp,clip]{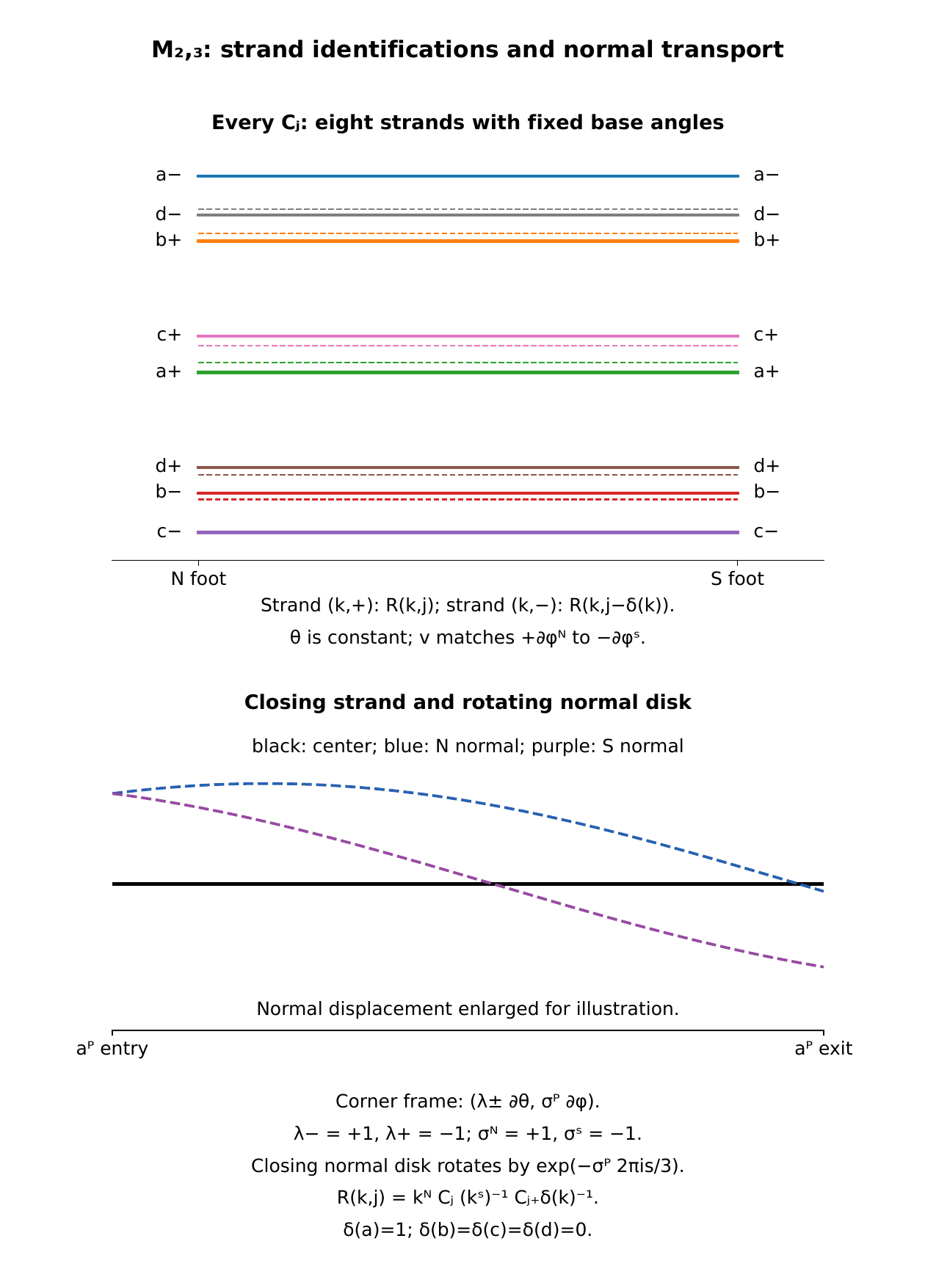}
\caption{Strand identifications through the connecting 1-handles
and normal transport for the closing 2-handle. The normal displacement
is enlarged in the drawing; the push-off uses a small displacement
with the same rotation. Product framings are specified by these
normal coordinates; see the framing convention at the start of Section~2.
Colors label local strands, not the global component colors of
Figure~\ref{fig:corners}.}
\label{fig:usualframematchings}
\end{figure}

Figure~\ref{fig:corners} projects the two distinct 0-handle boundaries.
The transparent foot outlines are graphical silhouettes, not certified
enclosures; strands may project across feet to which they do not attach.
Pair $k^-$ with $k^+$ within each pole and $C_j$ across poles, retaining
all eleven 1-handles. The twelve product curves are colored; black
intervals belong to the closing curves. Figures~\ref{fig:usualmatchings}
and~\ref{fig:usualframematchings} specify the strand identifications. In the base ribbon, a side at
$\theta_{k-}-\zeta$ matches $\theta_{k+}+\zeta$, and the other
side has the opposite signs. This is the side reversal that gives the
closing cycle $0,3,2,1,4,7,6,5$. Colors in Figures~\ref{fig:usualmatchings}
and~\ref{fig:usualframematchings} distinguish local strands;
they do not identify components in the full projections.

For the projection, regard $q\in S^3\subset\mathbb R^4$ as a
row vector. Choose a unit pole $\ell$ off the curves, a $4\times3$
matrix $\mathbf E$ whose columns are an orthonormal basis of
$\ell^\perp$, and a $3\times3$ orthogonal matrix $\mathbf O$.
The projected spatial coordinates are
\[
             X(q)=\frac{q\mathbf E}{1-q\cdot\ell}\,\mathbf O.
\]
The first two coordinates give the planar projection; the third
determines the overpasses. All curves in a drawing use this same map.
For these projections use $\epsilon=0.032$, $\zeta=0.012$, and
$\eta_{\rm proj}=0.005<\zeta$. These are smaller product neighborhoods
of the same framed cores. To compare with the preceding atlas, vary
the parameters through
\[
 (\epsilon_r,\zeta_r,\eta_r)
 =(1-r)(0.18,0.05,0.022)+r(0.032,0.012,0.005),
 \qquad 0\leq r\leq1.
\]
The foot angles remain fixed. Throughout this deformation,
$0<\epsilon_r<\pi/4$ and $0<\eta_r<\zeta_r$, and the foot neighborhoods
remain disjoint. Changing the cutoff extends each corner arc along
its fixed pair of arguments; changing the foot width extends the
closing intervals into the same collars. Reparameterize the adjacent
handle collars accordingly and carry the normal frames with them.
The product characteristic maps of the subsequent attaching spheres
are transported by the same collar isotopy. Hence this change preserves
the framed attaching maps and all higher attachments.
The camera data are recorded in
\path{figures/usual_projection_data.json}; the interval certificates are
in \path{certified_projection_crossings.json}.
The stored
camera entries are rounded decimals; orthonormalizing them changes
each entry by less than $10^{-12}$ and preserves all the crossings
and overpasses below. The attaching maps and framings are still
those of the product charts.

\begin{proposition}\label{prop:certifiedcrossings}
The analytic arcs in the six pole projections have only transverse
double points, with the following numbers of crossings:
\[
\begin{array}{c|rr}
 &N&S\\\hline
M_{2,3}&23&38\\
M_{3,2}&9&29\\
M_{2,4}&58&66
\end{array}
\]
Their overpasses agree with the recorded figure data. No individual
arc has a self-intersection and no three arcs have a common projected
point.
\end{proposition}
Appendix~\ref{app:crossings} proves Proposition~\ref{prop:certifiedcrossings}
by interval subdivision over the full parameter domains, including
endpoints. Numerical searches only choose candidate neighborhoods.


For $R_{k,j}$, let $q=j+\delta(k)$, with indices modulo three.
The attaching curve follows
\[
\begin{aligned}
 &\Gamma_{N,k-,q}^{-1},\quad k^N:q\longrightarrow j,\quad
 \Gamma_{N,k+,j},\quad C_j:N\longrightarrow S,\\
 &\Gamma_{S,k+,j}^{-1},\quad (k^S)^{-1}:j\longrightarrow q,\quad
 \Gamma_{S,k-,q},\quad C_q:S\longrightarrow N.
\end{aligned}
\]
At each connecting 1-handle $C_j$, strand $(k,+)$ belongs to $R_{k,j}$ and
strand $(k,-)$ belongs to $R_{k,j-\delta(k)}$. The strand's base angle
is constant from $N$ to $S$. These eight strands through each connecting 1-handle match the
same corner arcs at both ends.

\begin{figure}[!htbp]
\centering
\includegraphics[width=.98\textwidth]{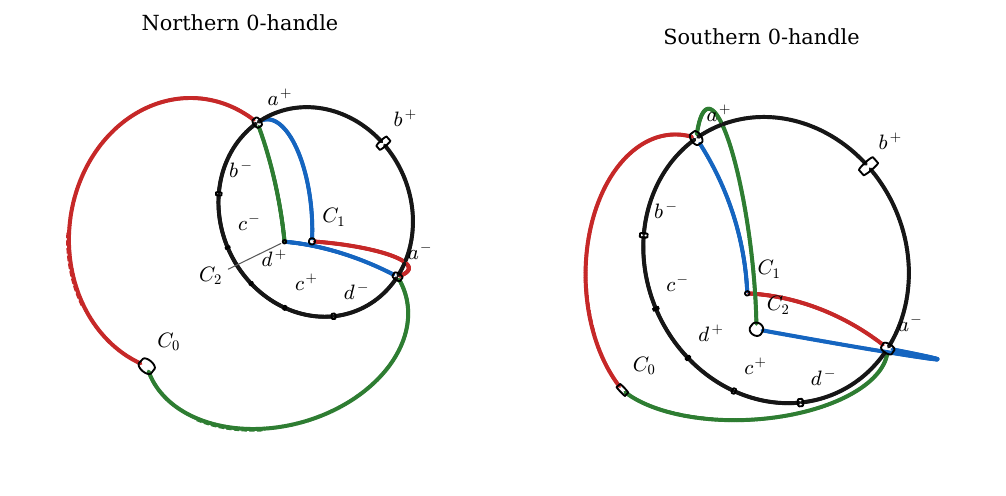}
\caption{The three rotating product 2-handles, shown as a sublink of
the preceding drawings. Black arcs belong to the closing 2-handles;
dashed arcs indicate the fiber-band normal direction, enlarged for visibility. Normal transport
is as in Figure~\ref{fig:usualframematchings}, with opposite
pole-coordinate signs.}
\label{fig:twelve}
\end{figure}

The remaining attaching curves are the boundaries of the two
base disk-and-four-band ribbons at pole coordinate zero. Their cycle
visits foot indices $0,3,2,1,4,7,6,5$, giving $[a,b][c,d]$.
Both sides of every band are drawn. The pole-coordinate normal disk
fixes the framing.
After the product 2-handles have been attached, the closing 2-handles
are $D^2_{\rm closing}\times D^2_N$ and
$D^2_{\rm closing}\times D^2_S$, with the pole disks as normal factors.
The parallel curves show the framing. Before the product 2-handles are attached, the normal framing
extending over each punctured pole section agrees with the Seifert
longitude of its boundary in the 0- and 1-handlebody. Here the closing
curve is null-homologous and the 0- and 1-handlebody has $H_2=0$.
Capping with the closing core gives the pole section, whose normal
Euler number is zero, since the commutator rotation paths cancel.
Each closing attachment therefore has coefficient zero relative to
this longitude, carried forward in a small neighborhood disjoint
from the other attachments. To obtain an integer label in a
dotted-circle presentation, this longitude must be carried through
the chosen conversion of the paired feet.
We may shrink the attaching
neighborhoods while transporting the product normal coordinates
and the face maps of the subsequent attaching spheres.  Along a base band, the
closing attachment follows the strips along the two boundary edges of the band, whereas
the product 2-handle attachments use strips in its interior.  On the base
0-handle these strips are separated by the corresponding foot-angle
intervals.  Taking disjoint product collars therefore separates the
closing attaching ribbons from the attaching ribbons of the product 2-handles.  The normal disk identifications for the higher attachments are taken
in the boundary after attaching the product 2-handles. They need not
extend as large normal disks in the original 0- and 1-handlebody.
Write $D_{P,j}=D^2_{\rm closing}\times\{q_{P,j}\}$ for its
disk face at pole $P$, where $q_{P,j}$ is the fiber-band foot
position in $D^2_P$. The two disk faces of the $j$th attaching
sphere are $D_{N,j}$ and $D_{S,j}$. Their collars are those of the
product construction.

\begin{theorem}\label{thm:attachingmaps}
The attaching maps above give the handle decomposition of $M_{2,3}$
in Theorem~\ref{thm:main}. The same construction gives the decompositions
of $M_{3,2}$ and $M_{2,4}$ drawn in the appendix.
\end{theorem}
\begin{proof}
The corner ribbons are disjoint by their two argument coordinates.
The three strands within each 1-handle $k^P$ remain at distinct disk angles,
and the eight strands in each connecting 1-handle have distinct base angles. The
specified endpoint matchings therefore give twelve disjoint framed
product 2-handle attaching curves. They are exactly the products of the four base bands
with the three fiber bands. The pole-center closing ribbons are
separate from those curves and attach the two products of the closing base disk with the pole disks. Finally the three spheres are
$\partial(D^2_{\rm closing}\times I_j)$, with disjoint $I_j$ collars
and product normal intervals. These are the remaining three product
attachments. Rounding the common faces recovers the quotient bundle.
The drawings show these attachments in separate boundary charts.
The scripts \path{check_attachment_atlases.py} and
\path{check_product_frames.py} check the finite attachment recipes
and symbolic frame identities. These consistency checks supplement
the product-chart argument; they do not independently establish the
geometric embeddings or frame transport.
For the other two pieces replace the base foot list by
$(a_i^-,b_i^+,a_i^+,b_i^-)_{i=1}^{g}$ and use
$\phi_{N,j}=2\pi j/p+\alpha$, $\phi_{S,j}=-\phi_{N,j}$.
There are $4gp$ corner charts at each pole, $2gp$ product 2-handles,
$4g$ strands through each connecting 1-handle, two closing 2-handles,
and $p$ attaching spheres for the 3-handles.
The same argument proves disjointness and compatibility of the framings.
The closing foot cycle is the concatenation of
$(4i,4i+3,4i+2,4i+1)$ for $0\leq i<g$, and reads
$\prod_{i=1}^g[a_i,b_i]$. The appendix draws these specializations.
\end{proof}

\subsection{Integral homology}
We compute the handle chain complex from the attaching curves and
the product faces of the attaching spheres.
\begin{theorem}\label{thm:piecehomology}
For $g\geq2$ and $p\geq2$, the displayed decomposition of $M_{g,p}$ has
\[
\begin{aligned}
H_0(M_{g,p};\Z)&=\Z,\\
H_1(M_{g,p};\Z)&=\Z^{2g}\oplus\Z/p,\\
H_2(M_{g,p};\Z)&=\Z^{1+2(g-1)(p-1)}\oplus\Z/p,\\
H_3(M_{g,p};\Z)&=H_4(M_{g,p};\Z)=0.
\end{aligned}
\]
\end{theorem}
\begin{proof}
Orient each connecting core $c_j$ from $N$ to $S$, and orient $R_{k,j}$ by
the displayed attaching path. Use the handle labels for the corresponding integral chain generators;
in particular, $c_j$ is the generator associated to $C_j$.
The 3-handle $\mathsf T_j$ has core $D^2_{\rm closing}\times I_j$.
The chain groups, in degrees zero through three, have ranks
$2$, $4g+p$, $2+2gp$, and $p$. Their differentials are
\begin{equation}\label{eq:diagramdifferentials}
\begin{aligned}
\partial_1(c_j)&=S-N,&\partial_1(k^P)&=0,\\
\partial_2(R_{k,j})&=k^N-k^S+c_j-c_{j+\delta(k)},&
\partial_2(H_P)&=0,\\
\partial_3(\mathsf T_j)&=H_S-H_N+R_{b_1,j-1}-R_{b_1,j}.&&
\end{aligned}
\end{equation}
Here all fiber indices are modulo $p$. The first two lines are the
signed handle passages of the attaching curves. For the last line,
the cap faces contribute $H_S-H_N$. Trace the annular face along
$[a_1,b_1]\cdots[a_g,b_g]$. Starting at fiber index $j$, positive
$a_1$ transport takes it to $j-1$. The first four faces contribute
\[
R_{a_1,j-1}+R_{b_1,j-1}-R_{a_1,j-1}-R_{b_1,j}.
\]
The $a_1$ faces cancel. All later commutators have trivial fiber
transport and their faces cancel in pairs. This proves the last line
from the attaching sphere itself. The formulas satisfy $\partial_1\partial_2=\partial_2\partial_3=0$.

A basis of $\ker\partial_1$ consists of the $4g$ base-handle
generators $k^N,k^S$ and the $p-1$ differences $c_j-c_0$, $1\leq j<p$. For $k\ne a_1$ the product 2-handle relations identify
$k^N$ and $k^S$. Put $\Delta_a=a_1^N-a_1^S$. The rotating 2-handle relations give
$c_{j+1}-c_j=\Delta_a$, and their cyclic sum gives $p\Delta_a=0$.
Eliminating $c_j-c_0=j\Delta_a$ and $a_1^S=a_1^N-\Delta_a$ leaves the $2g$ northern
base edges free and $\Delta_a$ of exact order $p$.

For degree two, a free integral basis of $\ker\partial_2$ is
\[
H_N,\ H_S,\quad z_{k,j}=R_{k,j}-R_{k,0}
\quad(k\ne a_1,\ 1\leq j<p).
\]
Indeed the connecting-handle coefficients force all coefficients of the rotating
$a_1$ 2-handles to be equal; the $a_1^N$ coefficient is then $p$ times
that integer and forces it to vanish. For every other $k$, the sum of
its product 2-handle coefficients must vanish, giving exactly the differences
$z_{k,j}$. Put $h=H_S-H_N$ and $z_j=z_{b_1,j}$, with $z_0=0$.
The 3-handle relations are
\[
                 h+z_{j-1}-z_j=0\qquad(j\bmod p).
\]
The relations give $z_j=jh$. Eliminating $z_1,\ldots,z_{p-1}$
leaves $ph=0$. The assignment $h\mapsto1$, $z_j\mapsto j$
to $\Z/p$ respects every relation, so $h$ has exact order $p$. The generator $H_N$ and all $z_{k,j}$ for
$k\notin\{a_1,b_1\}$ remain free. This gives the asserted $H_2$.
Finally a vanishing combination of the 3-handle columns forces all
its coefficients equal, by the cyclic $b_1$ differences. Its cap
coefficient is $p$ times that common integer, so it is zero.
Thus $\partial_3$ is injective and $H_3=0$. There are no 4-handles.
\end{proof}

The complement exact sequence gives an independent check of the
$H_2$ torsion; see Appendix~\ref{app:complementcheck}.

\begin{corollary}\label{cor:poletorsion}
Let $\mathcal S_N$ and $\mathcal S_S$ be the two pole sections of
$M_{g,p}\to\Sigma_g$, oriented by the base. Then
$[\mathcal S_S]-[\mathcal S_N]$ generates the torsion subgroup of
$H_2(M_{g,p};\Z)$ and has exact order $p$.
\end{corollary}
\begin{proof}
In the product handle complex the two sections have cellular
2-cycles $H_N$ and $H_S$. Thus their difference is the class $h$
in the preceding proof. Summing the $p$ attaching-sphere relations
cancels the $R_{b_1,j}$ terms and gives $ph=0$. Conversely,
eliminating $z_1,\ldots,z_{p-1}$ leaves exactly this relation on
$h$, so its order is $p$. This also identifies the torsion class
geometrically rather than only through a Smith normal form.
\end{proof}

\begin{center}
\begin{tabular}{ccc}
\toprule Piece & $H_1(-;\Z)$ & $H_2(-;\Z)$\\\midrule
$M_{2,3}$ & $\Z^4\oplus\Z/3$ & $\Z^5\oplus\Z/3$\\
$M_{3,2}$ & $\Z^6\oplus\Z/2$ & $\Z^5\oplus\Z/2$\\
$M_{2,4}$ & $\Z^4\oplus\Z/4$ & $\Z^7\oplus\Z/4$\\\bottomrule
\end{tabular}
\end{center}
The Smith-form checks in \path{check_diagram_homology.py}, recorded
in \path{verified_diagram_homology.json}, verify the chain algebra;
Theorem~\ref{thm:attachingmaps} supplies its geometric identification.

\section{Boundary markings and surface gluings}\label{sec:markings}
\subsection{The covering surface and boundary inclusion}\label{subsec:boundaryinclusion}
The boundary is connected even though the fiber has $p$ boundary
circles: $r_p$ permutes them transitively. It is a circle bundle over $F_{g,p}$, and its map to the base
factors through the covering $q_p:F_{g,p}\to\Sigma_g$.
The chosen normal framing gives
\[
 \partial M_{g,p}\cong F_{g,p}\times S^1\cong\Sigma_G\times S^1.
\]
Cut the base along a nonseparating curve dual to the character $\chi_p$,
take $p$ copies of the cut surface, and cyclically identify the cut
boundaries. Choose the basepoint lift in sheet zero, with the base paths
from this cut-and-glue model. This gives the marked surface-group generators:
\[
 E=a_1^p,\quad B=b_1,\quad
 C_{i,r}=a_1^r a_{i+1}a_1^{-r},\quad
 D_{i,r}=a_1^r b_{i+1}a_1^{-r}.
\]
Here $1\leq i<g$, $0\leq r<p$. With the ordering by sheet $r$, they
satisfy
\[
 [E,B]\prod_{r=0}^{p-1}\prod_{i=1}^{g-1}[C_{i,r},D_{i,r}]=1.
\]
To check this directly, write $B_r=a_1^r b_1a_1^{-r}$ and
$K_r=\prod_i[C_{i,r},D_{i,r}]$. Lifting the base relation gives
$B_{r+1}=K_r^{-1}B_r$ for $r<p-1$; the last lift gives
$[E,B_0]K_0\cdots K_{p-1}=1$. 

For the surface relator we order the pairs by sheet $r$.
For homology, put $(e_0,f_0)=([E],[B])$ and
$(e_{i,r},f_{i,r})=([C_{i,r}],[D_{i,r}])$; after the initial pair,
use lexicographic order $(i,r)$. The pairs are represented by the
cover's initial torus and its lifted one-holed tori. Reordering
the based generators requires Hurwitz exchanges; the resulting
conjugations are recorded in the mixed-family marking of
Section~\ref{subsec:mixed}. Define
$Q_{p,g}=(q_p)_*:H_1(F_{g,p};\mathbb Z)\to H_1(\Sigma_g;\mathbb Z)$.
Its values on this basis are
\[
 Q_{p,g}(e_0)=p a_1,\quad Q_{p,g}(f_0)=b_1,\qquad
 Q_{p,g}(e_{i,r})=a_{i+1},\quad Q_{p,g}(f_{i,r})=b_{i+1}.
\]

\begin{figure}[!htbp]
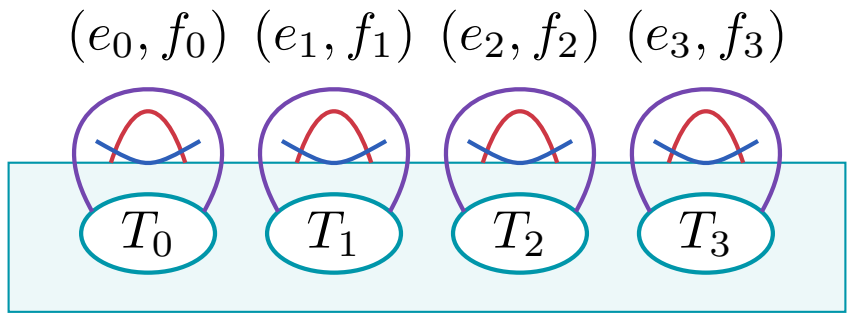

\centering
\surfacepairs{4}
\caption{A marked connected-sum model of the genus-four
boundary surface. Each $T_j$ is a one-holed torus attached to the
planar core, and carries its indicated red-blue symplectic marking. The drawing
records the surface marking used to define the gluing maps.
The torus outlines depict subsurfaces.}
\end{figure}

Choose the deleted disks to be geodesic disks of one fixed small
radius, and let $z_0$ be the center of the reference deleted disk,
whose boundary circle has based peripheral loop $x_1$. For $\eta\in S^1$, write
$\mathcal R_\eta=R_{\arg\eta}$. On the pullback to the
root cover $q_p:F_{g,p}\to\Sigma_g$, the map
\[
 F_{g,p}\times P_p\longrightarrow q_p^*M_{g,p},\qquad
 ((x,z),v)\longmapsto ((x,z),\mathcal R_{z/z_0}v)
\]
is a global fiberwise diffeomorphism. It fixes the pole section.
For the selected cuff, its boundary restriction identifies
$F_{g,p}\times S^1$ with $\partial M_{g,p}$, using the equatorial framing of~\cite{AB,AY}, with the oriented
equatorial tangent and transverse vector as the normal frame.

Fix a reference path $\gamma:[0,1]\to P_p$ from the north pole to
one point of the selected cuff. For a based root-cover loop
$\ell(s)=(x(s),z(s))$, transport the entire path by
$\mathcal R_{z(s)/z_0}$. The map
\[
 H(s,t)=\bigl(x(s),\mathcal R_{z(s)/z_0}\gamma(t)\bigr),
 \qquad 0\leq s,t\leq1,
\]
is a homotopy from the pole-section lift to the fixed-angle boundary
lift, with the two vertical sides equal to the chosen reference path.
In the pullback product trivialization, this path is constant in the
$P_p$ factor. Thus the boundary lift of $a_1^p$ has no additional
meridional factor, even when $z(s)$ winds around $S^1$. The same
argument applies to every covering generator, so their based images
are exactly the displayed words in the kernel of $\chi_p$.
Under the splitting of Theorem~\ref{thm:piecehomology}, the boundary inclusion induces
\begin{equation}\label{eq:boundaryinclusion}
\begin{aligned}
j_* &: H_1(F_{g,p};\Z)\oplus\Z
\longrightarrow H_1(\Sigma_g;\Z)\oplus\Z/p,\\
j_*(v,k)&=(Q_{p,g}v,k\bmod p).
\end{aligned}
\end{equation}
The cyclic target summand is generated by the positive meridian
$t=x_1$. Here positive means the boundary orientation of the oriented
normal disk to the multisection, opposite to the cuff orientation
induced from $P_p$. With the chosen connecting-handle orientations,
$[t]=[c_0-c_1]=-\Delta_a$.

The same trivialization proves boundary incompressibility.
The selected boundary component of the pullback has inclusion
\[
 \pi_1(F_{g,p})\times\mathbb Z
 \longrightarrow \pi_1(F_{g,p})\times\pi_1(P_p),
\]
equal to the identity in the first factor and the peripheral-circle
inclusion in the second. The latter is injective for $p\geq2$.
The pullback has $p$ boundary components, one for each reference cuff.
Each maps with degree one to the original boundary: a boundary point
determines its unique nearby root center, and the choice of reference
cuff then fixes the root-cover coordinate. Thus the selected component
maps diffeomorphically to the original boundary. The covering
homomorphisms are injective.
Thus the boundary inclusion is $\pi_1$-injective.

\subsection{Surface diffeomorphisms and symplectic compatibility}
Identify each cover with a surface assembled from a planar core and
ordered marked one-holed tori $T_0,\ldots,T_{G-1}$. This identification
uses the covering basis above.
Reindex the $G$ homology pairs as $(e_j,f_j)$, with the initial pair
indexed by zero and
\[
 j=1+(i-1)p+r
\]
for the pair previously indexed by $(i,r)$. Choose based torus
generators $A_j,B_j$ with $[A_j]=e_j$ and $[B_j]=f_j$, including the
Hurwitz base-path corrections required by this order. They satisfy
\[
 \prod_{j=0}^{G-1}[A_j,B_j]=1.
\]
The indices of adjacent tori below refer to this ordered based marking.
The mixed-family correction is given explicitly in
Section~\ref{subsec:mixed}.
Cap the outer boundary of the planar
core, keeping the basepoint in that cap. We use Dehn twists and
exchanges of adjacent one-holed tori.

For an integer $c$, on $T_j$ let $\mathcal U_{j,c}$ and $\mathcal L_{j,c}$ be the powers
of the twists, identity on the boundary collar, with actions
\[
 U_c=\begin{pmatrix}1&c\\0&1\end{pmatrix},\qquad
 L_c=\begin{pmatrix}1&0\\c&1\end{pmatrix}
\]
on the column basis $(e_j,f_j)$. The twist signs are fixed by these
matrices. On standard based torus generators $A,B$ one can choose
representatives
\[
 \mathcal U_c:(A,B)\mapsto(A,BA^c),\qquad
 \mathcal L_c:(A,B)\mapsto(AB^c,B).
\]
Both preserve $[A,B]$ exactly. All matrix products act on column
vectors, so the rightmost twist is performed first.

Define $\mathcal E_j$ to exchange the adjacent marked tori $T_j$ and
$T_{j+1}$. Move their attaching circles along complementary semicircles in the
planar core, preserving the chosen collar coordinates. Choose the signed exchange
and transported base paths by the Hurwitz convention below. Transport the attached tori by
their marked parameterizations. The resulting diffeomorphism is the
identity outside the genus-two subsurface containing these tori and their joining region and sends each marked
pair to the corresponding pair in its new position.
For based generators $A_j,B_j$, put $h_j=[A_j,B_j]$ and choose the
exchange sign by the following base path convention:
\[
\begin{aligned}
\mathcal E_j(A_j)&=A_{j+1},&\mathcal E_j(B_j)&=B_{j+1},\\
\mathcal E_j(A_{j+1})&=h_{j+1}^{-1}A_jh_{j+1},&
\mathcal E_j(B_{j+1})&=h_{j+1}^{-1}B_jh_{j+1}.
\end{aligned}
\]
All other generators are fixed. The Hurwitz move is
$h_jh_{j+1}\mapsto h_{j+1}(h_{j+1}^{-1}h_jh_{j+1})$,
which preserves the ordered surface relator exactly. The transported
base path of the second torus follows the inverse boundary loop of the
first. This path accounts for the conjugations in the formula.
Together with the twist formulas, these substitutions determine the
induced homomorphisms on the boundary fundamental groups.

\begin{figure}[!htbp]
\centering
\begin{tikzpicture}[>=Stealth,font=\small,scale=0.9]
\draw[fill=myturq!7,draw=myturq,rounded corners] (-3,-1.5) rectangle (3,1.6);
\draw[fill=white,draw=myred,thick] (-1.3,0) circle (0.43);
\draw[fill=white,draw=myblue,thick] (1.3,0) circle (0.43);
\node at (-1.3,0) {$T_j$};\node at (1.3,0) {$T_{j+1}$};
\draw[myred,very thick,->] (-1.3,0.48) .. controls (-1.3,1.4) and (1.3,1.4) .. (1.3,0.48);
\draw[myblue,very thick,->] (1.3,-0.48) .. controls (1.3,-1.3) and (-1.3,-1.3) .. (-1.3,-0.48);
\node[above] at (0,1.62) {exchange $\mathcal E_j$};
\node[below] at (0,-1.55) {the attached tori travel with their red-blue markings};
\end{tikzpicture}
\caption{Unbased trajectories of an exchange in the planar core.
The Hurwitz formulas specify the exchange sign and the transported base paths.}
\end{figure}
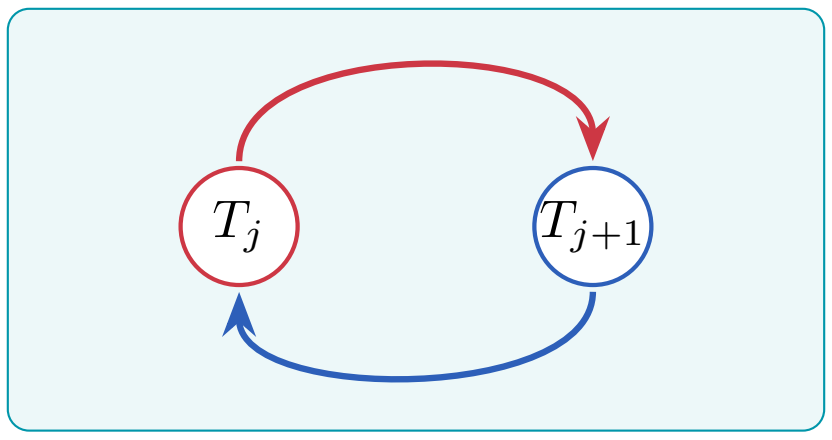

Let $\mathcal B_j$ be the $2\times2$ homology matrix assigned to the $j$th
torus. Its entries are specified for each family in Section~4.
Write $\pi(j)$ for the target position of this torus; a permutation
is listed as $(\pi(0),\ldots,\pi(G-1))$. Apply the specified maps on the marked tori
to obtain $D_m$, then perform the specified exchanges to obtain $P$. The subscript
$m\geq1$ is the integer twist parameter.
Our surface gluing is $f_m=P\circ D_m:F_L\to F_R$, where
$F_L,F_R$ are the two marked covering surfaces. Its action on first
homology is the $2G\times2G$ matrix $S_m=(f_m)_*$, with target/right
rows and source/left columns in the chosen ordered bases.
For $v\in\Sigma_G$ and normal-circle angle $\theta$, the boundary gluing is
\[
 \phi_m(v,e^{i\theta})=(f_m(v),e^{-i\theta}).
\]
It preserves the surface orientation and reverses the normal circle.
The normal-circle reversal is independent of $v$.

Rescale the product symplectic forms to give the multisections equal area.
Choose the marked identifications in local area coordinates so that
the induced forms agree near one fixed disk containing the cap,
foot patches, and tree; for the symmetric case use identical copies. Write $\omega_L,\omega_R$ for the induced forms in the chosen
markings. The tree joins the chosen foot patches and has a disk
neighborhood; it does not contain the essential cycles on the
one-holed tori. For Appendix~\ref{app:relative} it is the lifted tree
$a_0\cup a_1$; see the
\hyperref[app:fixedfeet]{fixed-disk placement in Appendix~\ref{app:relative}}.
The complement of this disk retains the genus. Choose the marked
tori and the paths of their exchanges in that complement.
The forms have equal total area and agree near the disk.
Choose all twists and exchanges to be the identity on this disk. Thus
$\beta=f_m^*\omega_R-\omega_L$ vanishes near the disk and has integral
zero. On the connected complement of a smaller disk still containing
the foot patches and tree, integration identifies compactly supported
degree-two cohomology with $\mathbb R$. Thus $\beta$, which has integral
zero, has a compactly supported primitive there. Extend this primitive
by zero to obtain $\alpha_{\rm M}$ with $d\alpha_{\rm M}=\beta$. The positive forms
$\omega_t=\omega_L+t\beta$, $0\leq t\leq1$, and the equation
$\iota_{X_t}\omega_t=-\alpha_{\rm M}$ give a flow $h_t$ fixed near these marked sets,
with $h_1^*f_m^*\omega_R=\omega_L$. Replacing $f_m$ by $f_m\circ h_1$
therefore preserves its based mapping class and foot charts while
making the surface gluing area preserving. Extend the isotopy by
\[
 (w,z)\longmapsto(h_t(w),z).
\]
The normal coordinate is fixed, so this extension preserves its
trivialization and introduces no surface-dependent normal-circle twist.
The vertical oriented planes project isomorphically to the normal
quotient of the multisection. Choose the transverse vector so that
the equatorial tangent and transverse vector form a positive frame.
Project this frame to the symplectic normal bundle and normalize by
positive rescaling, without an additional circle-valued change of frame.
Choose symplectic tubular charts in
this framing. Normal-circle reversal then corresponds to a complex
anti-linear identification of the trivial normal bundles, as required
in the symplectic sum construction~\cite{Gompf}.
The inverse of the adjusted surface map is $h_1^{-1}\circ f_m^{-1}$.
Transporting the entire attaching graph together with its framing through
$h_t^{-1}$ gives an isotopy fixed near the feet. Its product lift,
with the normal-circle coordinate fixed, extends over a boundary collar
by isotopy extension. These are the same framed attachments, relative
to the feet, as those
specified using $f_m^{-1}$ in Appendix~\ref{app:relative}.
Composing $f_m$ with a Torelli mapping class preserves $S_m$
but may change the smooth sum.

\subsection{The integral gluing calculation}
Let $M_L=M_{g,p}$ and $M_R=M_{h,q}$ with
$1+p(g-1)=1+q(h-1)=G$, and let
$\phi(v,e^{i\theta})=(f(v),e^{-i\theta})$ for an
orientation-preserving surface diffeomorphism $f$.
Using the equatorial framings, pole-section splittings, and inclusion
\eqref{eq:boundaryinclusion}, we recall the product-framed
Mayer--Vietoris formula to fix the boundary-map conventions.

\begin{proposition}[\cite{AB}, Proposition~3.7]\label{prop:generalHone}
Put $S=f_*$ and
\[
 A_f=\begin{pmatrix}Q_{p,g}\\-Q_{q,h}S\end{pmatrix}.
\]
For $Z=M_L\cup_\phi M_R$ there is an isomorphism
\[
 H_1(Z;\Z)\cong\coker A_f\oplus\Z/\gcd(p,q).
\]
Here $\Z/1=0$. The matrix $A_f$ has $2(g+h)$ rows and $2G$ columns,
so it is square precisely when $G=g+h$. If it is also nonsingular, then
\[
 |H_1(Z;\Z)|=|\det A_f|\,\gcd(p,q).
\]
\end{proposition}
\begin{proof}
The pieces and their common boundary are connected. The map on
$H_0$ in Mayer--Vietoris sends $1$ to $(1,-1)$ and is injective.
Consequently $H_1(Z;\Z)$ is the cokernel of the boundary map on
$H_1$. In the specified splittings this map is
\[
 (v,k)\longmapsto
 \bigl(Q_{p,g}v,-Q_{q,h}Sv;\ k\bmod p,k\bmod q\bigr).
\]
Indeed the minus sign in Mayer--Vietoris cancels the reversal of
the normal circle. The surface and meridian variables lie in
separate summands. The surface cokernel is $\coker A_f$, while
\[
 ((\Z/p)\oplus(\Z/q))/\langle(1,1)\rangle
 \cong\Z/\gcd(p,q).
\]
This proves the direct-sum formula and, in the square nonsingular
case, the order formula.
\end{proof}

A gluing of the form
$(v,e^{i\theta})\mapsto(f(v),\lambda(v)e^{-i\theta})$, with
$\lambda:\Sigma_G\to S^1$, can add a right-meridian term to the
surface columns through $\lambda_*$. The direct-sum calculation
above therefore uses the stated product normal-circle gluing; it
is not a consequence of the surface matrix alone.

\section{The two gluing families}\label{sec:families}\subsection{The symmetric \texorpdfstring{$(3,3)$}{(3,3)} case}
Both pieces are $M_{2,3}$, with base generators $(a,b,c,d)$ and
handle counts $(2,11,14,3)$. Only the $a$ face pairing rotates the
pair-of-pants fiber.

Put $x=x_1$, $y=x_2$ for the first two peripheral loops.
In these coordinates a useful presentation is
\[
\begin{split}
\pi_1(M_{2,3})=\langle a,b,c,d,x,y\mid{}&[a,b][c,d]=1,\\
&axa^{-1}=y,\quad aya^{-1}=(xy)^{-1},\\
&[s,x]=[s,y]=1\quad(s=b,c,d)\rangle.
\end{split}
\]
The third peripheral loop is $z=(xy)^{-1}$. This presentation follows from the handle attachments.
In the disk-and-band basis the same monodromy is
$y_1\mapsto y_2y_1^{-1}$, $y_2\mapsto y_1^{-1}$.

Abbreviate $C_r=C_{1,r}$ and $D_r=D_{1,r}$.
The marked genus-four covering surface has pairs
\[
 (E,B),\quad(C_0,D_0),\quad(C_1,D_1),\quad(C_2,D_2),
\]
and surface relation $[E,B][C_0,D_0][C_1,D_1][C_2,D_2]=1$.
Thus $(A_0,B_0)=(E,B)$ and
$(A_{r+1},B_{r+1})=(C_r,D_r)$ for $0\leq r<3$.
With the based paths from the covering polygon, its boundary inclusion
is specified by the following words:
\begin{center}
\begin{tabular}{cll}
\toprule Pair & First generator & Second generator\\\midrule
$0$ & $E\mapsto a^3$ & $B\mapsto b$\\
$1$ & $C_0\mapsto c$ & $D_0\mapsto d$\\
$2$ & $C_1\mapsto aca^{-1}$ & $D_1\mapsto ada^{-1}$\\
$3$ & $C_2\mapsto a^2ca^{-2}$ & $D_2\mapsto a^2da^{-2}$\\\bottomrule
\end{tabular}
\end{center}
The meridian maps to $x$. All these words commute with $x$ because
their monodromy is trivial. Thus the inclusion respects the boundary relation
$[\pi_1(\Sigma_4),x]=1$.

In the homology order
$(e_0,f_0,e_1,f_1,e_2,f_2,e_3,f_3)$, the projection matrix is
\[
 Q_{3,2}=\begin{pmatrix}
3&0&0&0&0&0&0&0\\
0&1&0&0&0&0&0&0\\
0&0&1&0&1&0&1&0\\
0&0&0&1&0&1&0&1
\end{pmatrix}.
\]
Use the same construction and marking for the second piece.
To obtain a presentation of the fundamental group of the symplectic sum,
combine the two piece presentations, identify the eight surface boundary
generators by $f_m$, defined below, and identify the meridians inversely.

\par\medskip
\textbf{The chosen genus-four gluing.}
Use the following torus permutation and blocks from~\cite[Theorem~4.5]{AY}:
\begin{center}
\begin{tabular}{cccc}
\toprule Source $j$ & Target $\pi(j)$ & $\mathcal B_j$ & Twist prescription\\\midrule
0&3&$U_1$&$\mathcal U_{0,1}$\\
1&0&$U_1$&$\mathcal U_{1,1}$\\
2&2&$U_{-1}$&$\mathcal U_{2,-1}$\\
3&1&$L_1U_{m-1}$&$\mathcal U_{3,m-1}$, then $\mathcal L_{3,1}$\\\bottomrule
\end{tabular}
\end{center}
For the permutation perform exchanges with indices
\[
 0,\ 2,\ 1,\ 2\qquad\text{in this order.}
\]
The ordered torus labels change as
\[
\begin{split}
(0,1,2,3)&\longrightarrow(1,0,2,3)\longrightarrow(1,0,3,2)\\
&\longrightarrow(1,3,0,2)\longrightarrow(1,3,2,0).
\end{split}
\]
The source-to-target permutation is $\pi=(3,0,2,1)$.
The based generator actions of the twists and exchanges in Section~3
therefore specify $f_m$ on the genus-four surface group, not only on
homology. For $m=1$, all eight images under $f_1^{-1}$ are displayed
in equation~\eqref{eq:inversegenusfour} of Appendix~\ref{app:relative}.

The $8\times8$ matrix $S_m$ has block $\mathcal B_j$ in row-pair
$\pi(j)$ and column-pair $j$, and zeros elsewhere. The surface
Mayer-Vietoris matrix, with left base rows followed by right base
rows and the boundary basis as columns, is
\[
A_m=\begin{pmatrix}Q_{3,2}\\-Q_{3,2}S_m\end{pmatrix}
=\begin{pmatrix}
3&0&0&0&0&0&0&0\\
0&1&0&0&0&0&0&0\\
0&0&1&0&1&0&1&0\\
0&0&0&1&0&1&0&1\\
0&0&-3&-3&0&0&0&0\\
0&0&0&-1&0&0&0&0\\
-1&-1&0&0&-1&1&-1&1-m\\
0&-1&0&0&0&-1&-1&-m
\end{pmatrix}.
\]
Its determinant is $-9m$. For the specified boundary gluing, the meridian relations give
an additional direct summand $\Z/3$.

\par\medskip
\begin{proposition}\label{prop:symmetric}
For every $m\geq1$,
$H_1(Z_m^{3,3};\Z)\cong(\Z/3)^2\oplus\Z/(3m)$.
\end{proposition}
\begin{proof}
Let $a,b,c,d$ denote the left base homology generators and $A,B,C,D$
the right generators. The columns of $A_m$ give the relations
\[
\begin{array}{ll}
3a-C=0,& b-C-D=0,\\
c-3A=0,&d-3A-B=0,\\
c-C=0,&d+C-D=0,\\
c-C-D=0,&d+(1-m)C-mD=0.
\end{array}
\]
Subtracting the fifth relation from the seventh gives $D=0$.
Consequently
\[
c=C=3A=3a,\qquad d=-C,\qquad b=C,\qquad B=-6A.
\]
The last relation gives $mC=0$. Eliminating generators by these coefficient-one relations leaves
\[
 3(a-A)=0,\qquad 3mA=0.
\]
Let $t,s$ be the left and right meridian classes. They satisfy $3t=3s=0$ and $t+s=0$; hence one independent
meridian class of order three remains. Therefore
\[
 H_1(Z^{3,3}_m;\Z)\cong(\Z/3)^2\oplus\Z/(3m).
\]
This agrees with~\cite{AY} for the surface gluing chosen here.
For $m=1,2$, we obtain
\[
 H_1(Z^{3,3}_1;\Z)=(\Z/3)^3,\qquad
 H_1(Z^{3,3}_2;\Z)=(\Z/3)^2\oplus\Z/6.
\]
\end{proof}

Use the left complement's absolute handles, insert a collar
$\Sigma_4\times S^1\times[0,1]$ with transition $\phi_m$, and attach
the right complement upside down relative to its entire boundary.
Its index-$k$ handles become relative index-$(4-k)$ handles.

The unreduced closed counts are $(2,14,28,14,2)$, so $\chi=4$.
For $m=1,2$, the first-homology orders $27$ and $54$ distinguish
the diffeomorphism types. The signature and rational cohomology
ring are computed in Section~5.

\par\medskip
For $m=1$, the right complement also admits a relative product
presentation with three 1-handles, fourteen 2-handles, eleven 3-handles,
and two 4-handles; its relative attachment charts are in
Appendix~\ref{app:relative}. This is an alternative to the upside-down
absolute presentation used for the counts above.

\subsection{The mixed \texorpdfstring{$(2,4)$}{(2,4)} case}\label{subsec:mixed}
The pieces $M_{3,2}$ and $M_{2,4}$ are respectively annulus and
four-holed-sphere bundles over $\Sigma_3$ and $\Sigma_2$; their common
boundary surface has genus five.

Let $\mu$ be the positively oriented peripheral loop with its chosen base path in the annulus;
it also generates its fundamental group. With base generators
$(a,b,c,d,e,f)$, the left group is
\[
\begin{split}
\pi_1(M_{3,2})=\langle a,b,c,d,e,f,\mu\mid{}&
[a,b][c,d][e,f]=1,\\
&a\mu a^{-1}=\mu^{-1},\quad[s,\mu]=1\ (s=b,c,d,e,f)\rangle.
\end{split}
\]
For the right base $(A,B,C,D)$ use $x_1,x_2,x_3,x_4$ with
$x_1x_2x_3x_4=1$. Its group is
\[
\begin{split}
\pi_1(M_{2,4})=\langle A,B,C,D,x_1,\ldots,x_4\mid{}&
[A,B][C,D]=1,\quad x_1x_2x_3x_4=1,\\
&Ax_jA^{-1}=x_{j+1},\quad[s,x_j]=1\ (s=B,C,D)\rangle,
\end{split}
\]
where $x_5=x_1$. Eliminating $x_4$ gives seven generators.
The two pieces have unreduced counts
\[
 (2,14,14,2),\qquad (2,12,18,4).
\]

\par\medskip
\textbf{The two genus-five boundary markings.}
For the left covering surface use the pair order
\[\begin{gathered}
0=(E,B_0),\quad1=(C_0,D_0),\quad2=(C_1,D_1),\\
3=(E_0,F_0),\quad4=(E_1,F_1).
\end{gathered}\]
Here $(E,B_0)$ is the initial covering pair; $(C_r,D_r)$ and
$(E_r,F_r)$ are the lifts of the second and third base pairs,
with the base-path correction given below. Put $W=[aca^{-1},ada^{-1}]$. The images in the left complement are:
\begin{center}
\begin{tabular}{cll}
\toprule Pair&First generator&Second generator\\\midrule
0&$a^2$&$b$\\
1&$c$&$d$\\
2&$aca^{-1}$&$ada^{-1}$\\
3&$W^{-1}eW$&$W^{-1}fW$\\
4&$aea^{-1}$&$afa^{-1}$\\\bottomrule
\end{tabular}
\end{center}
The original sheet-first closing product is
\[
 [a^2,b][c,d][e,f]W[aea^{-1},afa^{-1}].
\]
In the displayed order its middle factors are instead
$[c,d]W(W^{-1}[e,f]W)$, which equal $[c,d][e,f]W$.
Hence the surface relation holds in this order. The base path of pair three
winds along $W^{-1}$ before reaching its original torus; omitting
this winding changes the based boundary map even though $Q_{2,3}$
is unchanged.

The right marking has initial pair $(A^4,B)$, followed by
$(A^rCA^{-r},A^rDA^{-r})$ for $r=0,1,2,3$. The left meridian is the based, positively oriented annulus
peripheral loop $\mu=x_1$, with $[\mu]=-\Delta_a$,
and the right meridian is $x_1$.
The covering polygon specifies the surface marking and its base paths.

In the homology order $(e_0,f_0,\ldots,e_4,f_4)$ the projection
matrices are
\[
Q_{2,3}=\begin{pmatrix}
2&0&0&0&0&0&0&0&0&0\\
0&1&0&0&0&0&0&0&0&0\\
0&0&1&0&1&0&0&0&0&0\\
0&0&0&1&0&1&0&0&0&0\\
0&0&0&0&0&0&1&0&1&0\\
0&0&0&0&0&0&0&1&0&1
\end{pmatrix},
\]
\[
Q_{4,2}=\begin{pmatrix}
4&0&0&0&0&0&0&0&0&0\\
0&1&0&0&0&0&0&0&0&0\\
0&0&1&0&1&0&1&0&1&0\\
0&0&0&1&0&1&0&1&0&1
\end{pmatrix}.
\]

\par\medskip
\textbf{The chosen genus-five gluing.}
Use the following source blocks and destinations from \cite{AY}:
\begin{center}
\begin{tabular}{cccc}
\toprule Source $j$ & Target $\pi(j)$ & $\mathcal B_j$ & Twist prescription\\\midrule
0&3&$L_1$&$\mathcal L_{0,1}$\\
1&1&$U_{-1}L_{m-1}$&$\mathcal L_{1,m-1}$, then $\mathcal U_{1,-1}$\\
2&2&$L_{-1}$&$\mathcal L_{2,-1}$\\
3&4&$L_{-1}$&$\mathcal L_{3,-1}$\\
4&0&$U_{-1}$&$\mathcal U_{4,-1}$\\\bottomrule
\end{tabular}
\end{center}
Perform adjacent exchanges with indices
\[
3,\ 2,\ 1,\ 0,\ 1,\ 2\qquad\text{in this order.}
\]
The final order is $(4,1,2,0,3)$. Thus $\pi=(3,1,2,4,0)$ gives
the source-to-target permutation.

In this genus-five case, $S_m$ is the $10\times10$ matrix
with block $\mathcal B_j$ in row-pair $\pi(j)$ and column-pair $j$.
The surface presentation matrix is
\[
\begin{split}
A_m&=\begin{pmatrix}Q_{2,3}\\-Q_{4,2}S_m\end{pmatrix}\\
&=\left(\begin{array}{rrrrrrrrrr}
2&0&0&0&0&0&0&0&0&0\\
0&1&0&0&0&0&0&0&0&0\\
0&0&1&0&1&0&0&0&0&0\\
0&0&0&1&0&1&0&0&0&0\\
0&0&0&0&0&0&1&0&1&0\\
0&0&0&0&0&0&0&1&0&1\\
0&0&0&0&0&0&0&0&-4&4\\
0&0&0&0&0&0&0&0&0&-1\\
-1&0&m-2&1&-1&0&-1&0&0&0\\
-1&-1&1-m&-1&1&-1&1&-1&0&0
\end{array}\right).
\end{split}
\]
Its determinant is $8m$. For this gluing, the meridian relations give an additional direct
summand $\Z/2$. Both the sign of
the normal circle and the pair convention are the same as in the
symmetric case.

\par\medskip
\begin{proposition}\label{prop:mixed}
For every $m\geq1$,
$H_1(Z_m^{2,4};\Z)\cong(\Z/2)^2\oplus\Z/(4m)$.
\end{proposition}
\begin{proof}
Let $a,b,c,d,e,f$ be the left base homology generators and $A,B,C,D$
the right ones. The ten columns of the preceding matrix give
\[
\begin{array}{ll}
2a-C-D=0,&b-D=0,\\
c+(m-2)C+(1-m)D=0,&d+C-D=0,\\
c-C+D=0,&d-D=0,\\
e-C+D=0,&f-D=0,\\
e-4A=0,&f+4A-B=0.
\end{array}
\]
The fourth minus the sixth relation gives $C=0$. Then
\[
b=d=f=D,\qquad c=e=-D,\qquad D=-4A,\qquad B=0.
\]
The third relation gives $mD=0$ and the first gives $2a=D$. Put $w=a+2A$;
the remaining surface relations are precisely
\[
 2w=0,\qquad4mA=0.
\]
For the left and right meridian classes $t,s$, the relations
$2t=0$, $4s=0$, $t+s=0$ leave one independent
order-two class. Hence
\[
 H_1(Z^{2,4}_m;\Z)\cong(\Z/2)^2\oplus\Z/(4m),
\]
and specifically
\[
 H_1(Z^{2,4}_1;\Z)=(\Z/2)^2\oplus\Z/4,\qquad
 H_1(Z^{2,4}_2;\Z)=(\Z/2)^2\oplus\Z/8.
\]
\end{proof}
Use the same collar construction with $\Sigma_5$ in place of
$\Sigma_4$, $M_{3,2}$ on the left, and $M_{2,4}$ on the right. Its
unreduced closed handle counts are
\[
 (h_0,h_1,h_2,h_3,h_4)=(2,18,32,14,2),
 \qquad2-18+32-14+2=4.
\]
For $m=1,2$, the first-homology orders are $16$ and $32$.

\section{Topology of the symplectic sums}
Write $Z_m$ for either family $Z_m^{3,3}$ or $Z_m^{2,4}$.
\begin{proposition}\label{prop:aspherical}
Let $M_L=M_{g,p}$ and $M_R=M_{h,q}$, where $g,h,p,q\geq2$
and their boundary surfaces have the same genus. Gluing these pieces
by any orientation-reversing boundary diffeomorphism gives an
aspherical closed manifold with infinite fundamental group.
In particular, this holds for both families $Z_m^{3,3}$ and
$Z_m^{2,4}$, for every $m\geq1$.
\end{proposition}
\begin{proof}
The cyclic cover of each complement associated to $\ker\chi_p$ is
$\Sigma_G\times P_p$, where $G=1+p(g-1)\geq2$.
Both factors are aspherical. Each boundary component of this cover is
$\Sigma_G\times S^1$, and its fundamental group injects because a
peripheral circle of $P_p$ has infinite order. Boundary injectivity
descends to $M_{g,p}$, recovering the incompressibility result
proved in~\cite{AB}. The chosen boundary diffeomorphism
gives an amalgamated product with injective edge maps,
\[
 \pi_1(Z_m)=\pi_1(M_L)*_{\pi_1(\Sigma_G)\times\Z}\pi_1(M_R).
\]

The universal cover is a tree of universal covers of the pieces and
boundary collars, over the Bass-Serre tree of the amalgam. Its vertex and edge spaces are contractible.
The lifted boundary collars make the edge inclusions cofibrations.
The union over a finite subtree is contractible by induction: adjoining
a terminal vertex gives a pushout of contractible spaces along a
contractible cofibration, hence a homotopy pushout with the homotopy
type of a point.

Choose a finite triangulation of the closed manifold compatible with
the splitting collar. Its lifted triangulation is locally finite, so
the image of any sphere map is contained in a finite subcomplex and
meets only finitely many lifted pieces and collars. Adjoin the finitely
many paths joining their vertices in the tree. The sphere map then
lies in the contractible union over that finite subtree and is
null-homotopic. Thus the universal cover is contractible by Whitehead's
theorem, and $Z_m$ is aspherical. The amalgam normal-form theorem
embeds the infinite boundary group in $\pi_1(Z_m)$.
\end{proof}

\begin{corollary}\label{cor:minimal}
Every symplectic sum in Theorem~\ref{thm:main} is symplectically minimal.
\end{corollary}
\begin{proof}
In an aspherical manifold every map from $S^2$ is null-homotopic,
so every embedded sphere represents zero in integral homology and
has self-intersection zero. In particular there is no embedded
sphere of square $-1$.
\end{proof}

Since $\chi(\partial M_L)=0$, additivity gives
$\chi(Z_m)=\chi(M_L)+\chi(M_R)=4$.
The removed tubular neighborhood is $\Sigma_G\times D^2$ and has
signature zero. Novikov additivity applied to the ambient product ruled
surface gives $\sigma(M_{g,p})=\sigma(\Sigma_g\times S^2)=0$.
A second application across the gluing boundary gives $\sigma(Z_m)=0$.

\begin{corollary}\label{cor:homology}
For either family and every $m\geq1$, put $T=H_1(Z_m;\Z)$.
Then
\[
\begin{aligned}
H_0(Z_m;\Z)&=H_4(Z_m;\Z)=\Z, & H_1(Z_m;\Z)&=T,\\
H_2(Z_m;\Z)&\cong\Z^2\oplus T, & H_3(Z_m;\Z)&=0.
\end{aligned}
\]
The rational cohomology ring is that of $S^2\times S^2$.
\end{corollary}
\begin{proof}
Since $T=H_1(Z_m;\Z)$ is finite, integral Poincare duality and the
universal coefficient theorem give
\[
 H_3(Z_m;\Z)\cong H^1(Z_m;\Z)=\operatorname{Hom}(T,\Z)=0.
\]
The Euler characteristic is four, so the free rank of $H_2$ is two.
The universal coefficient theorem gives
$\operatorname{Tor}H^2(Z_m;\Z)=\operatorname{Ext}(T,\Z)\cong T$,
noncanonically. Poincare duality therefore gives
$H_2(Z_m;\Z)\cong\Z^2\oplus T$. Connectedness and orientation give
$H_0=H_4=\Z$.

The free intersection pairing is unimodular of rank two and signature
zero, hence has determinant $-1$. It is hyperbolic over $\Q$:
for a matrix $\left(\begin{smallmatrix}a&b\\b&c\end{smallmatrix}\right)$,
the equality $b^2-ac=1$ gives a nonzero rational isotropic vector $u$.
Choose $v$ with $u\cdot v=1$, and replace $v$ by
$v-\tfrac12(v\cdot v)u$. The resulting basis has
matrix $\left(\begin{smallmatrix}0&1\\1&0\end{smallmatrix}\right)$.
Together with $H^1(Z_m;\Q)=H^3(Z_m;\Q)=0$, this identifies the graded
rational cohomology ring with that of $S^2\times S^2$.
\end{proof}
We do not determine whether these sums admit complex structures.

Appendix~\ref{app:relative} describes the relative product attachments
for the symmetric sum. The smooth sums obtained from other Torelli
lifts have not been compared.

\FloatBarrier
\appendix
\section{Handlebodies of \texorpdfstring{$M_{3,2}$ and $M_{2,4}$}{M(3,2) and M(2,4)}}\label{app:othercomplements}
\subsection{The complement \texorpdfstring{$M_{3,2}$}{M(3,2)}}
Use the genus-three base with six bands and two fiber bands interchanged
by the monodromy. In the figures,
$(a_1,b_1,a_2,b_2,a_3,b_3)=(a,b,c,d,e,f)$.
There are twenty-four corner charts at each pole, twelve product
2-handles, two closing 2-handles, and two 3-handles.
The chart conventions and normal frames are those of Section 2.
\begin{figure}[!htbp]
\centering
\includegraphics[width=\textwidth]{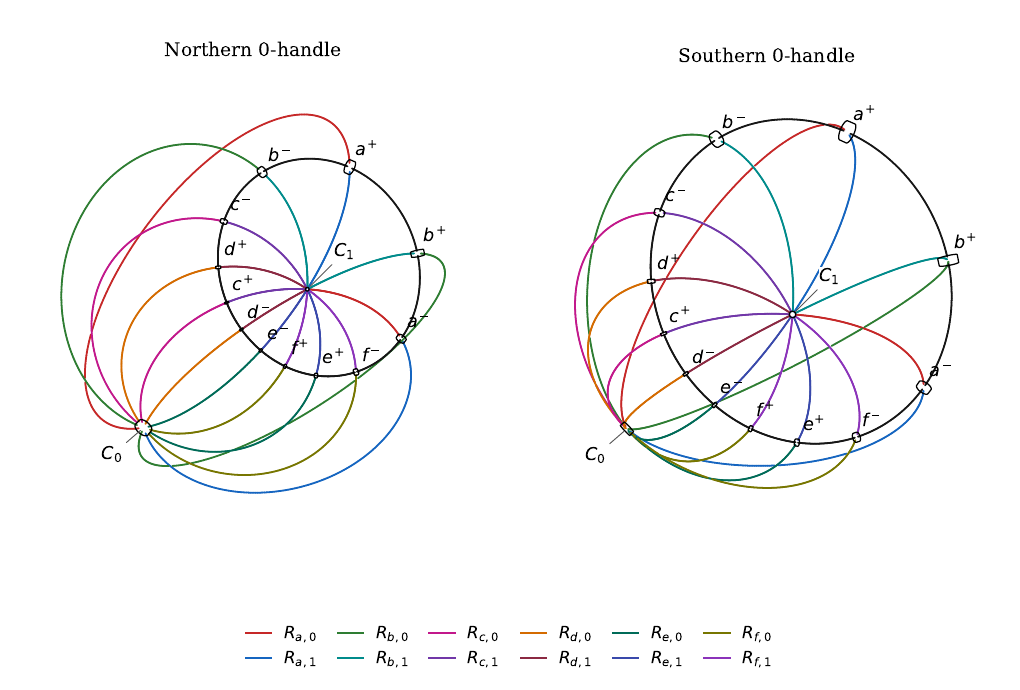}
\caption{$M_{3,2}$. Northern and southern projections of the product
and closing 2-handle attaching curves. Pair $k^-$ with $k^+$ at each
pole and $C_j$ across poles. Black arcs belong to the closing
2-handles. Figure~\ref{fig:M32transport} specifies the strand identifications and product framings.}
\end{figure}
\begin{figure}[!htbp]
\centering
\includegraphics[width=.84\textwidth,height=.36\textheight,keepaspectratio]{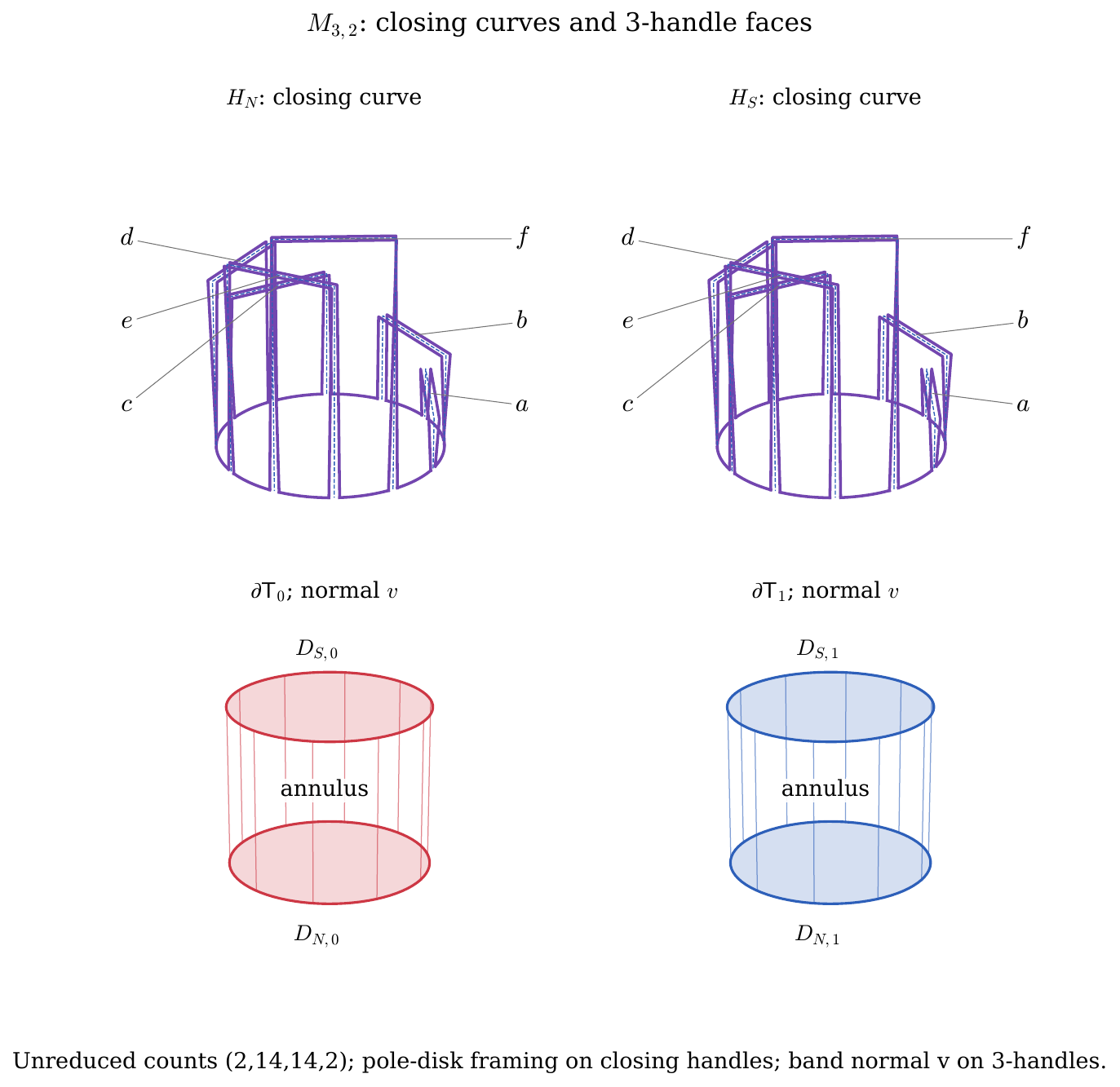}
\caption{$M_{3,2}$. The two closing 2-handle attaching curves and the product faces of the subsequent 3-handle attaching spheres, with their normal coordinates. These are local product schematics; their colors distinguish the displayed faces.}
\end{figure}
\begin{figure}[!htbp]
\centering
\includegraphics[width=\textwidth,height=.82\textheight,keepaspectratio,trim=4.00bp 4.66bp 5.90bp 4.00bp,clip]{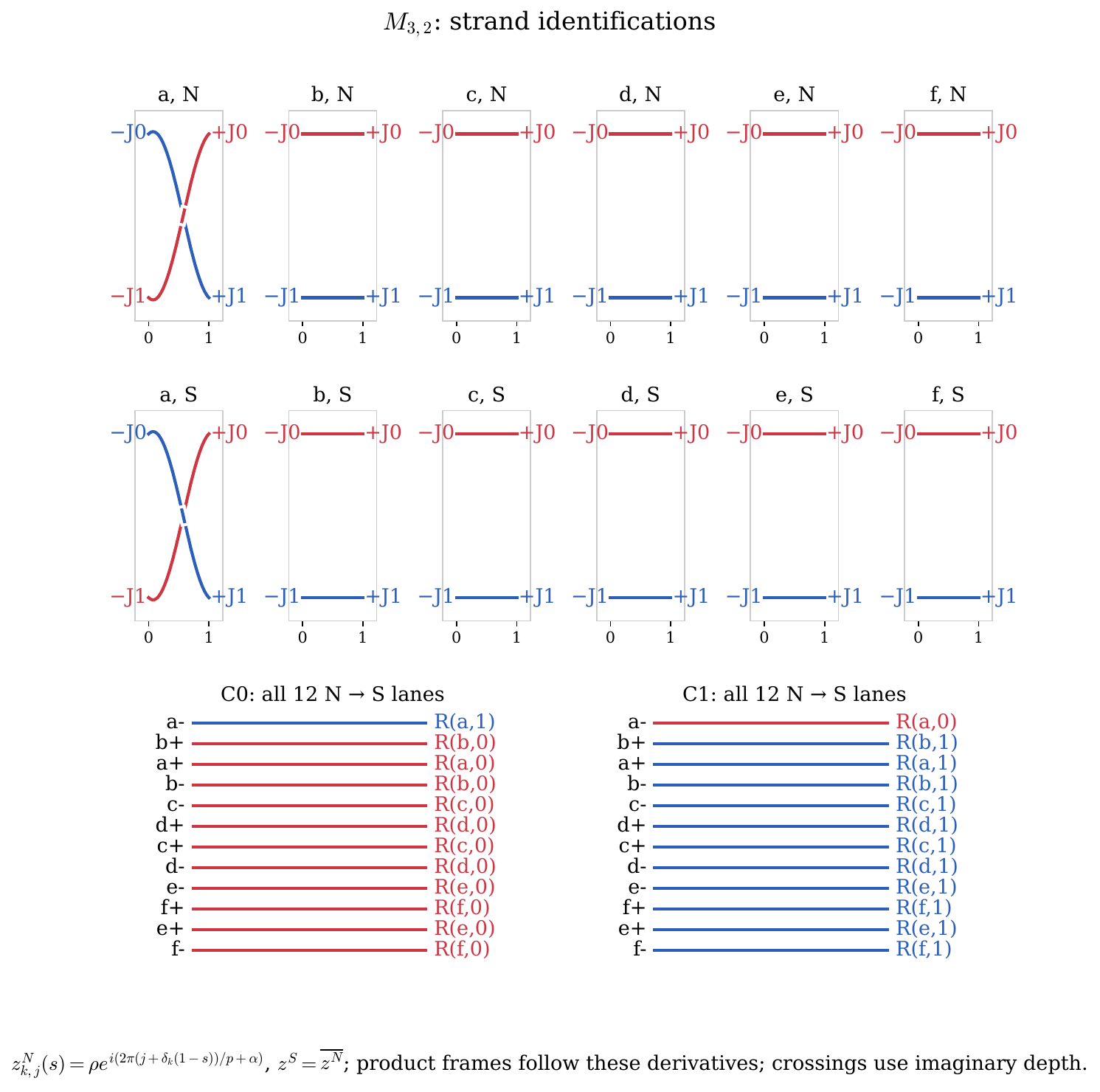}
\caption{$M_{3,2}$. Strand identifications through the 1-handles, with the normal coordinates inducing their framings.}
\label{fig:M32transport}
\end{figure}

\subsection{The complement \texorpdfstring{$M_{2,4}$}{M(2,4)}}
Use the genus-two base with four bands and four fiber bands permuted cyclically by the monodromy.
All figures use $(a,b,c,d)$ for the ordered base generators
$(a_1,b_1,a_2,b_2)$; these correspond to $(A,B,C,D)$ in Section~4.
There are thirty-two corner charts at each pole, sixteen product
2-handles, two closing 2-handles, and four 3-handles.
\begin{figure}[!htbp]
\centering
\includegraphics[width=\textwidth]{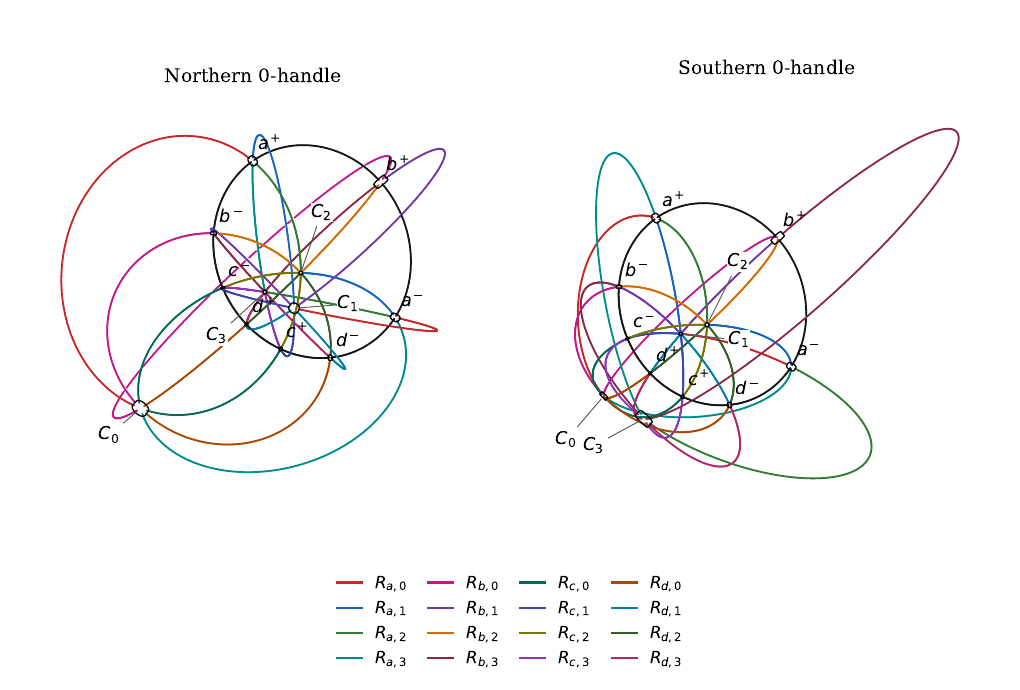}
\caption{$M_{2,4}$. Northern and southern projections of the product
and closing 2-handle attaching curves. Pair $k^-$ with $k^+$ at each
pole and $C_j$ across poles. Black arcs belong to the closing
2-handles. Figure~\ref{fig:M24transport} specifies the strand identifications and product framings.}
\end{figure}
\begin{figure}[!htbp]
\centering
\includegraphics[width=.84\textwidth,height=.36\textheight,keepaspectratio]{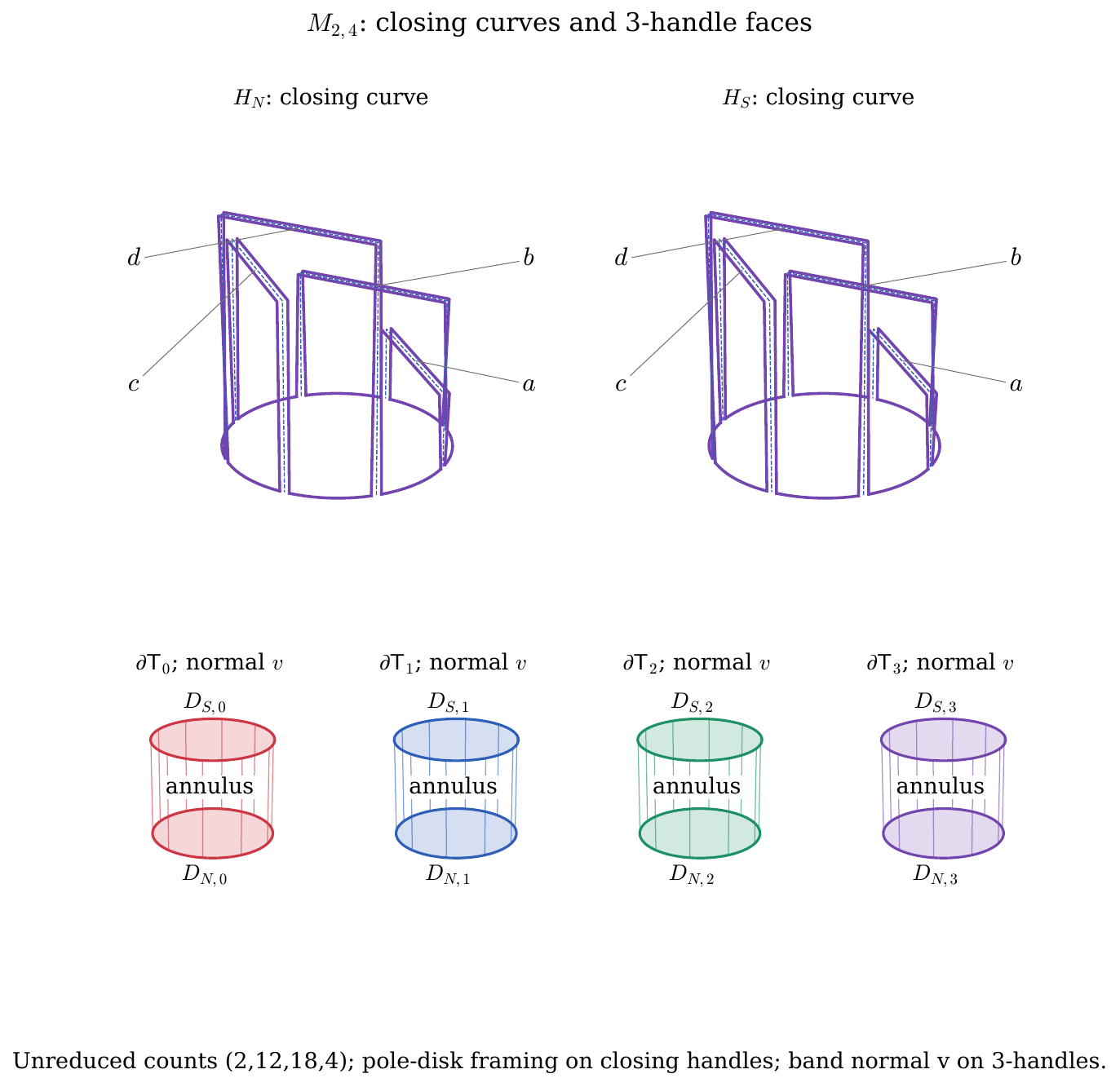}
\caption{$M_{2,4}$. The two closing 2-handle attaching curves and the product faces of the subsequent 3-handle attaching spheres, with their normal coordinates. These are local product schematics; their colors distinguish the displayed faces.}
\end{figure}
\begin{figure}[!htbp]
\centering
\includegraphics[width=\textwidth,height=.82\textheight,keepaspectratio,trim=4.00bp 4.66bp 5.90bp 4.00bp,clip]{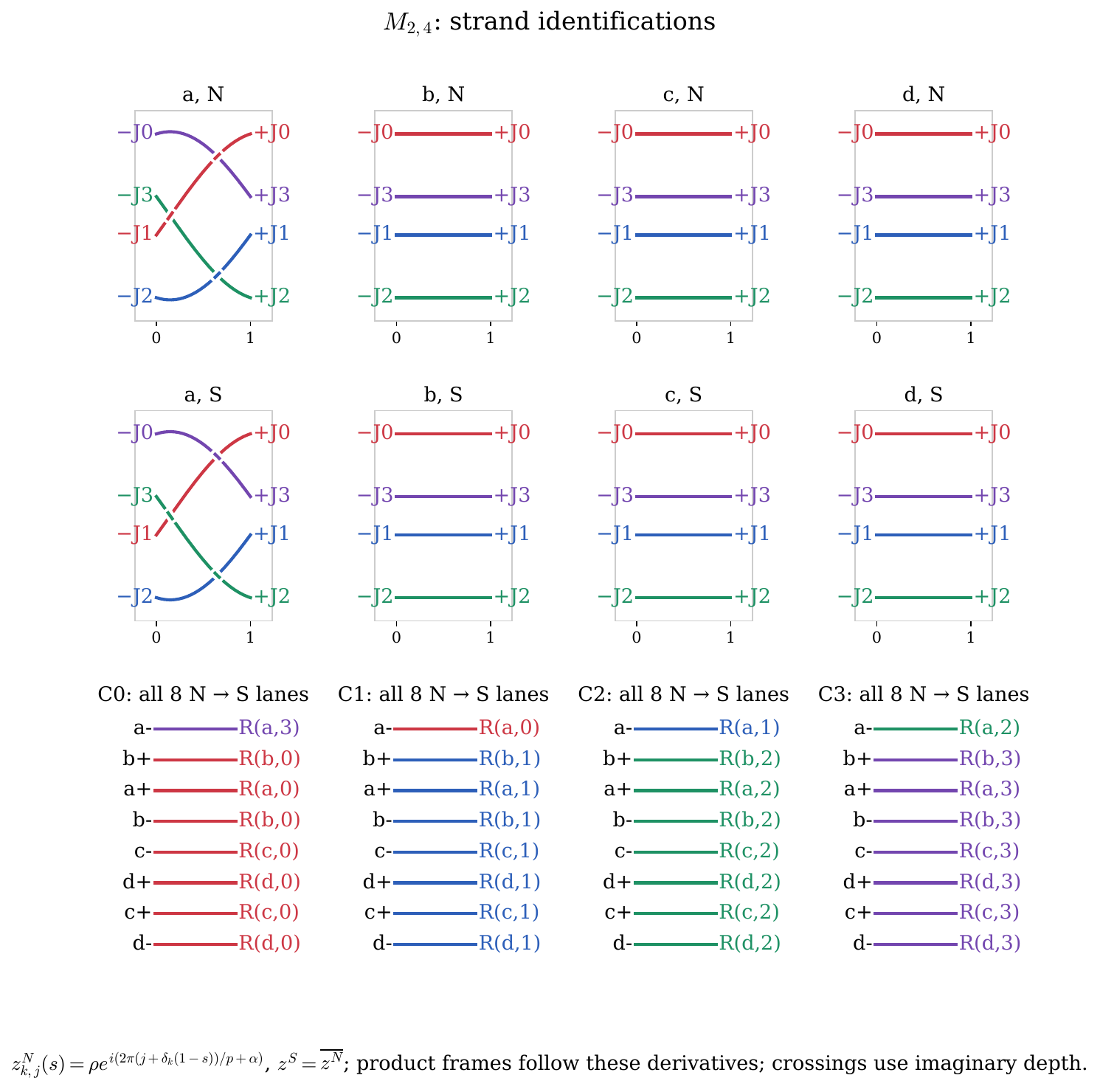}
\caption{$M_{2,4}$. Strand identifications through the 1-handles, with the normal coordinates inducing their framings.}
\label{fig:M24transport}
\end{figure}

\section{Relative attachments for the symmetric sum}\label{app:relative}
We describe the right-hand $M_{2,3}$ relative to its entire boundary
and pull its attachments back by $\phi_1^{-1}$ to form $Z_1^{3,3}$.
Unlike this relative product construction, the dual decomposition
recovers the original 2-handle belt circles only after adding the
three 1-handles dual to the original 3-handles; those circles need
not survive intact in the final boundary.

\textbf{The relative fiber.}
Use longitude $\theta$ and latitude $\lambda$ on the unit sphere.
Choose a longitude origin $\alpha_0$ ($\alpha_0=\pi/3$ in the
figure) and put
$\alpha_j=\alpha_0+2\pi j/3$, with $\alpha_3=\alpha_0+2\pi$.
Take geodesic disks of radius $0<\varepsilon<\pi/6$ about the equatorial
points $z_j=(\cos\alpha_j,\sin\alpha_j,0)$ of the sphere. Its framed
boundary circles, parametrized by the meridional angle
$\psi\in\mathbb R/2\pi\mathbb Z$, are
\[
 w_j(\psi)=\cos\varepsilon\,z_j+
 \sin\varepsilon\,(\cos\psi\,\tau_j+\sin\psi\,n),
\]
where $\tau_j=(-\sin\alpha_j,\cos\alpha_j,0)$ and $n=(0,0,1)$.
Put $\delta=\varepsilon/2$ and
$\Delta(\lambda)=\arccos(\cos\varepsilon/\cos\lambda)$.
The invariant band between holes $j$ and $j+1$ is
\[
 c_j(u,v):\quad\lambda=\delta u,\qquad
 \theta=\frac{\alpha_j+\alpha_{j+1}}2+
 v\left(\frac\pi3-\Delta(\delta u)\right),
 \qquad -1\leq u,v\leq1.
\]
Increasing $u$ runs from south to north. The absolute band core,
oriented from north to south in Section~2, is therefore traversed in
decreasing $u$, with tangent $-\partial_u$. This constant reversal
introduces no winding of the product normal frame. The sides $v=-1$
and $v=1$ lie on cuffs $j$ and $j+1$, respectively: their equations
are $\cos\lambda\cos(\theta-\alpha_j)=\cos\varepsilon$ and
$\cos\lambda\cos(\theta-\alpha_{j+1})=\cos\varepsilon$.
The cocore arc of this surface 1-handle is
$\eta_j(v)=c_j(0,v)$, with endpoints
$w_j(0)$ and $w_{j+1}(\pi)$. The three cocore arcs cut the fiber into
north and south disks. In the relative decomposition, these cocore arcs serve as the cores
of three relative surface 1-handles $E_j$, followed by two relative surface 2-handles. Their attaching curves
$\Theta_N,\Theta_S$ use the two edges of each relative 1-handle and the complementary
upper and lower cuff arcs, as drawn below.

\begin{figure}[!htbp]
\centering
\includegraphics[width=.90\textwidth,trim=4.00bp 4.13bp 3.91bp 4.00bp,clip]{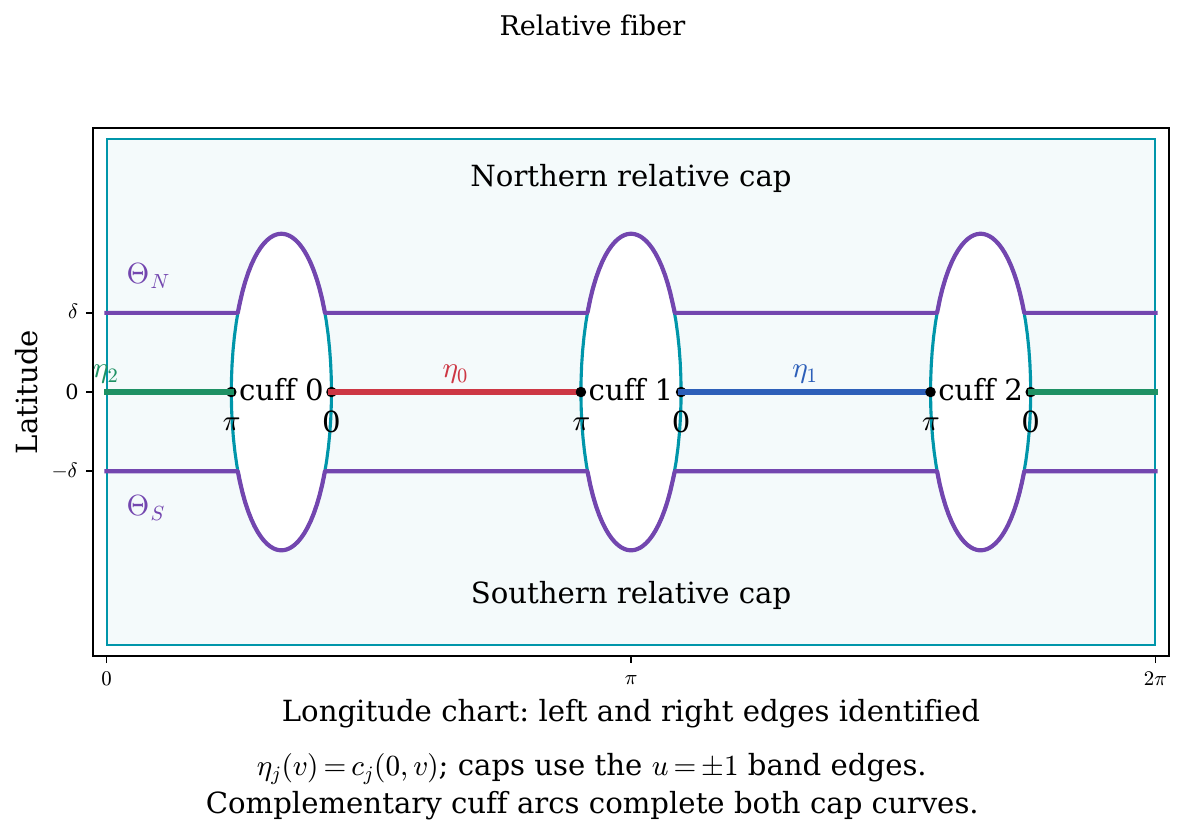}
\caption{The relative fiber in longitude and latitude coordinates,
with left and right edges paired. The colored arcs are cocores of the
absolute fiber bands and cores of the relative 1-handles. The purple
curves attach the two relative surface 2-handles. Their products with
the base 0-handle give four-dimensional 2-handles, with the base disk
as normal factor. The longitude origin is chosen so that all three
deleted disks lie inside the chart.}
\label{fig:relativecaps}
\end{figure}

In the higher-handle figures, $B_2=D^2_B=D^2_{\rm closing}$ is
the base closing disk, and $D^2_N,D^2_S$ denote the relative caps
with boundaries $\Theta_N,\Theta_S$. Their polar coordinates are
$(r_B,\theta_B)$ and $(r_{\rm cap},\theta_{\rm cap})$.
The intervals $I_s,I_v$ follow the indicated core coordinates;
$e_j$ in these figures is the core of $E_j$.
Multiplying by one base disk, four base bands, and a closing base disk
gives the relative handle counts
\[
 (h_0,h_1,h_2,h_3,h_4)=(0,3,14,11,2).
\]
The fourteen 2-handles consist of two products of the base 0-handle with the relative fiber caps and twelve products
of a base 1-handle with a relative surface 1-handle. The latter have core coordinates
$(s,v)$ and normal coordinates $(t,u)$. The eleven 3-handles are the
three products of the base 2-handle with relative surface 1-handles and eight base-band/cap products.
Their attaching spheres are the rounded boundaries of the product
cores; the remaining interval is the normal coordinate. The two
closing-disk/cap products are the final 4-handles.

As a check on these relative counts, Poincar\'e--Lefschetz duality
and Theorem~\ref{thm:piecehomology} give
\[
\begin{aligned}
H_4(M_{2,3},\partial M_{2,3};\Z)&=\Z,&
H_3(M_{2,3},\partial M_{2,3};\Z)&=\Z^4,\\
H_2(M_{2,3},\partial M_{2,3};\Z)&\cong\Z^5\oplus\Z/3,&
H_1(M_{2,3},\partial M_{2,3};\Z)&=\Z/3,
\end{aligned}
\]
with $H_0(M_{2,3},\partial M_{2,3};\Z)=0$.
Both alternating sums are two.

\begin{figure}[!htbp]
\centering
\includegraphics[width=\textwidth,height=.58\textheight,keepaspectratio,trim=3.00bp 4.70bp 2.75bp 2.00bp,clip]{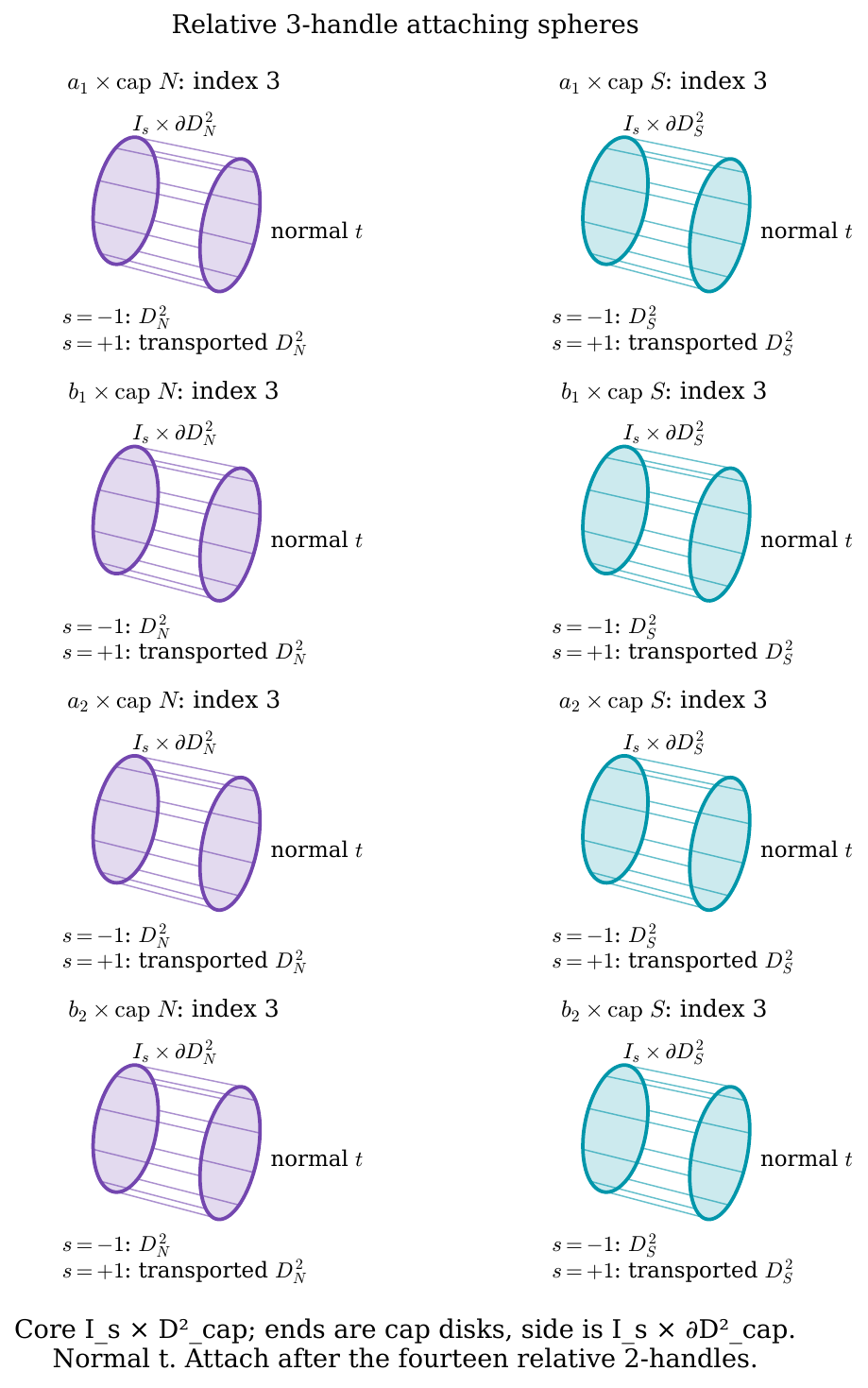}
\caption{The attaching spheres of the eight relative 3-handles obtained from
base bands and fiber caps. Each core is $D^1\times D^2$, with the base-band transverse interval as normal. The faces are identified according to the relative product decomposition. Labels $a_1,b_1,a_2,b_2$ denote $a,b,c,d$.}
\end{figure}
\begin{figure}[!htbp]
\centering
\includegraphics[width=.72\textwidth,trim=3.00bp 2.14bp 2.56bp 2.00bp,clip]{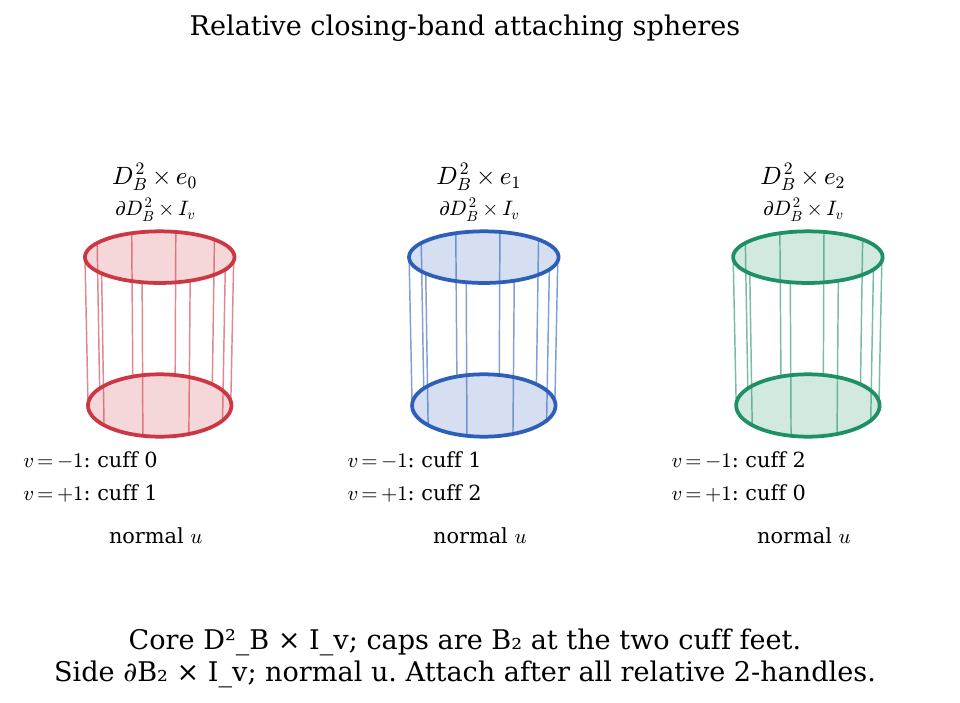}\par\vspace{5pt}
\includegraphics[width=.60\textwidth,trim=2.00bp 1.78bp 3.37bp 2.00bp,clip]{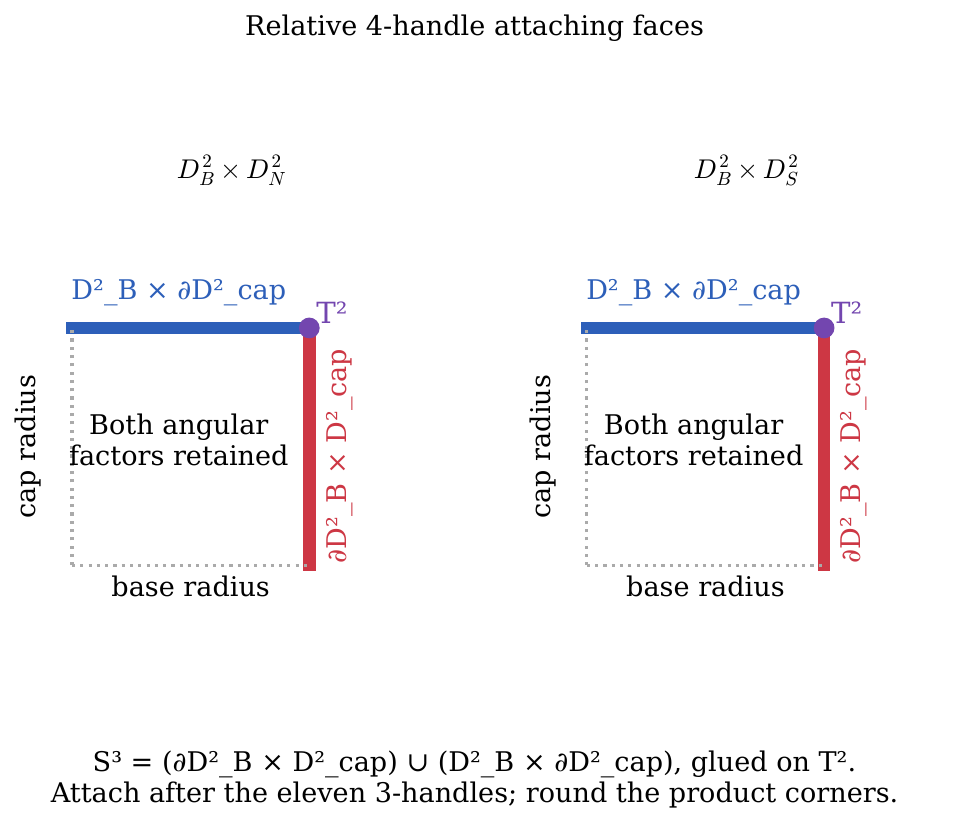}
\caption{Above: the attaching spheres for the three products of
the base surface 2-handle with relative surface 1-handles. Below:
the boundaries of the two 4-handle cores, each obtained by joining
two solid tori along their common boundary. The labels
$a_1,b_1,a_2,b_2$ denote $a,b,c,d$.}
\end{figure}

\phantomsection\label{app:fixedfeet}
\textbf{The surface marking and the six feet.}
Let $v_j$ be the three surface lifts above the center of the chosen
base disk. Place their small patches in one fixed disk of the planar
core of the marked genus-four surface. Include a regular neighborhood
of the lifted tree $a_0\cup a_1$ in this disk, carrying the covering
marking and its base paths by the same surface isotopy. The remaining,
essential parts of the lifted graph lie outside the disk.
Take the twists and exchanges
to be supported outside this disk. The resulting surface
diffeomorphism induces the gluing map on first homology specified in~\cite{AY}. The six relative 1-handle feet are
\[
 (v_j,\psi=0),\qquad(v_{j+1},\psi=\pi),\quad j\bmod3.
\]
Write $V_r$ for the vertex on sheet $r\in\mathbb Z/3$ of the lifted base polygon. The feet are $(V_{-j},0)$ and
$(V_{-(j+1)},\pi)$. The inverse gluing fixes their surface patches and reverses $\psi$.
Thus it fixes each foot center and reflects its latitude coordinate
in the foot chart $(x,y,u)$, where $x,y$ are local coordinates
on the covering surface and $u$ is the transverse fiber coordinate. Extend this reflection across each
relative 1-handle by
$(v,x,y,u)\mapsto(v,x,y,-u)$. Thus the core coordinate $v$ and
base-disk coordinates $(x,y)$ are fixed. The cap ribbons in
Figure~\ref{fig:relativecaps} are then interchanged by the circle
reversal; their base-disk normal frames are preserved.

\begin{figure}[!htbp]
\centering
\includegraphics[width=.68\textwidth,trim=4.50bp 3.89bp 5.34bp 4.00bp,clip]{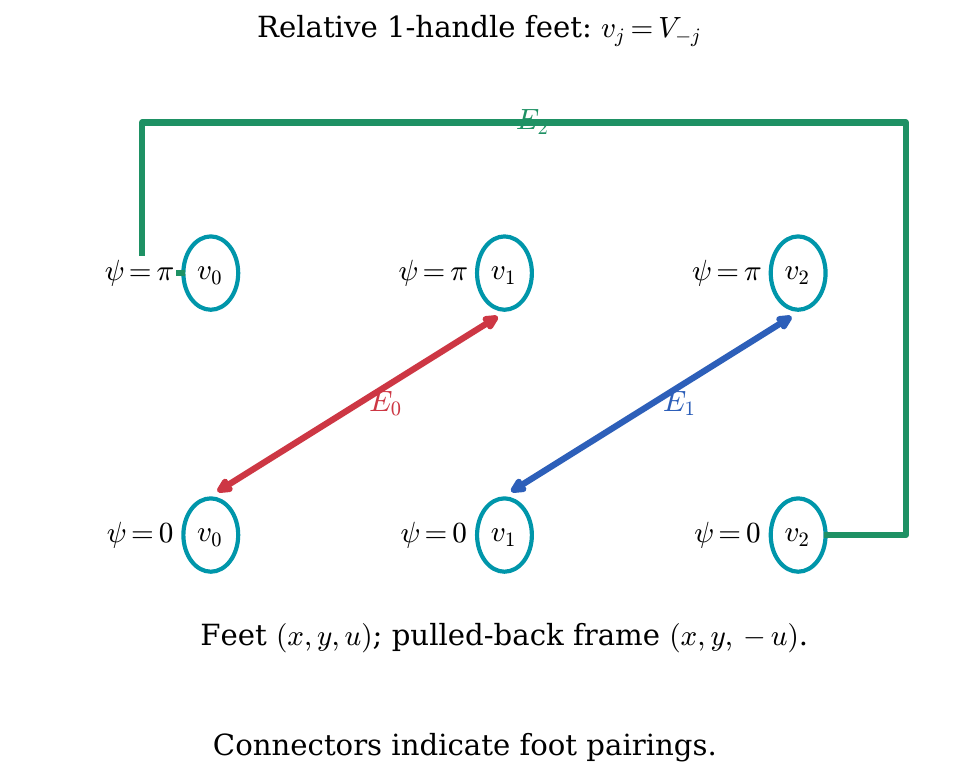}
\caption{The six marked feet of the three relative 1-handles for
the closed symmetric sum, with labels $v_j=V_{-j}$.
The connectors indicate the foot pairings; they are not 2-handle
attaching arcs. The displayed foot reflection and
fixed surface patches specify their attaching 3-ball parameterizations.}
\label{fig:relativefeet}
\end{figure}

\begin{figure}[!htbp]
\centering
\includegraphics[width=.95\textwidth,height=.54\textheight,keepaspectratio,trim=5.00bp 4.36bp 6.30bp 4.00bp,clip]{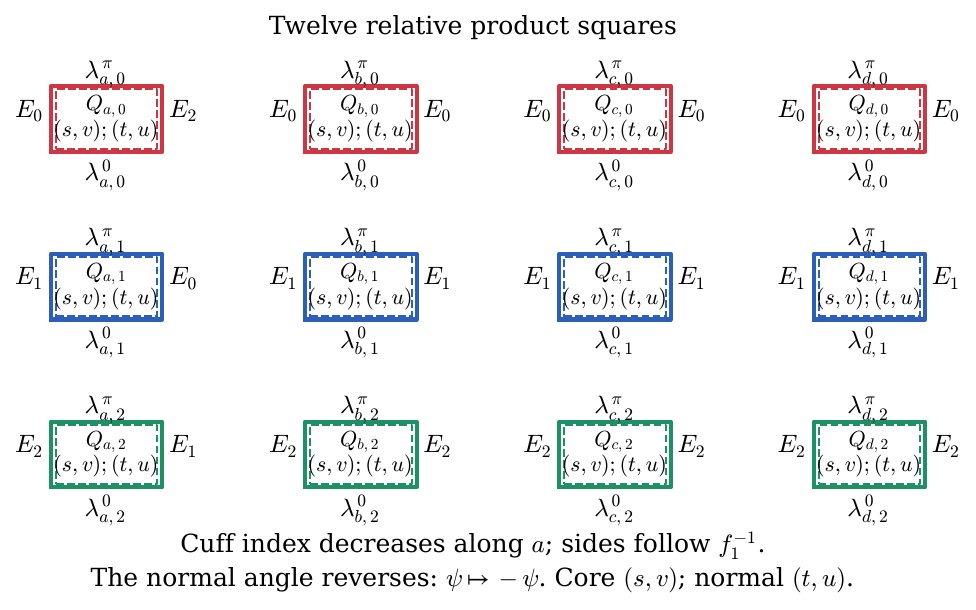}
\caption{The twelve relative product 2-handle charts, showing how the
sides attach to the relative 1-handles and the boundary surface. The sides on the boundary surface are transported by the inverse
gluing map specified below. These are local product charts. The
framing is induced by the transported normal coordinates $(t,u)$.}
\label{fig:relativesquares}
\end{figure}

\textbf{The twelve relative product 2-handles for $m=1$.}
Cuff indices decrease along the positively oriented $a$ 1-handle.
Use Schreier sheet index $r=-j\bmod3$ on the covering surface:
the positive $a$ band sends $j$ to $j-1$ and $r$ to $r+1$.
The cuff patches and Schreier vertices satisfy $v_j=V_{-j}$.
For $k\in\{a,b,c,d\}$, let $k_r$ be its lift starting at $V_r$.
In particular, $a_r$ runs from $V_r$ to $V_{r+1}$.
Choose $a_0,a_1$ as a tree in the fixed cap. The lifted graph, with its endpoints fixed, has the following
based edge labels:
\[
\begin{array}{c|ccc}
 &r=0&r=1&r=2\\\hline
a_r&1&1&A_0\\
b_r&B_0&h_1^{-1}B_0&h_2^{-1}h_1^{-1}B_0\\
c_r&A_1&A_2&A_3\\
d_r&B_1&B_2&B_3
\end{array}
\qquad h_i=[A_i,B_i].
\]
The entries $1$ for tree edges mean that their closed based labels
are trivial; the embedded tree edges and endpoints are retained.
Let $Q_{k,j}$ denote the relative product 2-handle obtained from
the base band $k$ and relative fiber handle $E_j$. Here
$\delta(a)=1$ and $\delta(b)=\delta(c)=\delta(d)=0$. Its two sides on the boundary surface are the lifts of
$k$ joining the relative 1-handle feet on sheets $-j$ and $-(j+1)$.
In Figure~\ref{fig:relativesquares}, these sides are labeled
$\lambda^0_{k,j}$ and $\lambda^\pi_{k,j}$, respectively, before
transport by the inverse surface gluing. Its sides in the relative 1-handles are $E_j$ and
$E_{j-\delta(k)}$. With endpoints fixed, the attaching path is
\[
 k_{-j}\,E_{j-\delta(k)}\,(k_{-(j+1)})^{-1}\,E_j^{-1}.
\]
These curves bound the product core squares, with the side curves
and their endpoints fixed by the product construction.

With $h_i=[A_i,B_i]$, the based substitutions for $T=f_1^{-1}$ are
\begin{equation}\label{eq:inversegenusfour}
\begin{aligned}
T(A_0)&=h_0A_1h_0^{-1},&T(B_0)&=h_0B_1A_1^{-1}h_0^{-1},\\
T(A_1)&=h_0h_2A_3B_3^{-1}h_2^{-1}h_0^{-1},&
T(B_1)&=h_0h_2B_3h_2^{-1}h_0^{-1},\\
T(A_2)&=h_0A_2h_0^{-1},&T(B_2)&=h_0B_2A_2h_0^{-1},\\
T(A_3)&=A_0,&T(B_3)&=B_0A_0^{-1}.
\end{aligned}
\end{equation}
The final pulled-back attaching path is
\[
 T(k_{-j})\,E_{j-\delta(k)}\,T(k_{-(j+1)})^{-1}\,E_j^{-1}.
\]
Here the $E$ labels are retained, with their foot parametrizations
reflected as above; $T$ acts on the marked surface paths.
Apply $T$ to the based edge labels in the table and reverse $\psi$ along
their cuff lifts. Apply the following maps, in order,
\[
 \mathcal E_2^{-1},\ \mathcal E_1^{-1},\ \mathcal E_2^{-1},\
 \mathcal E_0^{-1},\ \mathcal L_{3,-1},\
 \mathcal U_{2,1},\ \mathcal U_{1,-1},\ \mathcal U_{0,-1}
\]
to one embedded lifted ribbon graph. Transport the framing with this embedded graph.
The individual based words do not specify the framed attachments.

\textbf{A simultaneous surface diagram.}
Let $D$ be the fixed disk containing the tree and the three vertex
patches, and put $S=\Sigma_4\setminus\Int D$.
After contracting $a_0\cup a_1$ within $D$, number the twenty ports
in their positive boundary order on $D$ as follows:
\[
\begin{array}{c|rrrrrrrrrr}
i&0&1&2&3&4&5&6&7&8&9\\\hline
p_i&a_2^+&b_2^-&b_1^-&b_0^-&a_2^-&b_0^+&c_0^+&d_0^-&c_0^-&d_0^+\\[2pt]
i&10&11&12&13&14&15&16&17&18&19\\\hline
p_i&b_1^+&c_1^+&d_1^-&c_1^-&d_1^+&b_2^+&c_2^+&d_2^-&c_2^-&d_2^+
\end{array}
\]
Here $+$ and $-$ denote the initial and terminal ends of an oriented
edge. The order on $\partial S$ is the reverse order.
Deleting the $b_1,b_2$ ends leaves the four consecutive blocks
$(A_i^+,B_i^-,A_i^-,B_i^+)$, $0\leq i\leq3$.
Thus the eight selected edges give the standard generators for the genus-four surface.
Choose the boundary basepoint in the sector of the third lifted
polygon face, with port paths along the boundary of the tree
neighborhood avoiding that sector. Equivalently, extend the tree
neighborhood by a small tongue into that face.
The first two face relations are
$b_0^{-1}h_1b_1=1$ and $b_1^{-1}h_2b_2=1$. They give
$b_1=h_1^{-1}B_0$ and $b_2=h_2^{-1}h_1^{-1}B_0$
in the free group $\pi_1(S)$ with these port paths; the third records
$R=\prod_{i=0}^3[A_i,B_i]$ as the boundary word.
This specifies the punctured-surface marking used in the edge table,
including its boundary paths.

\begin{figure}[!htbp]
\centering
\includegraphics[width=.98\textwidth,height=.69\textheight,keepaspectratio]{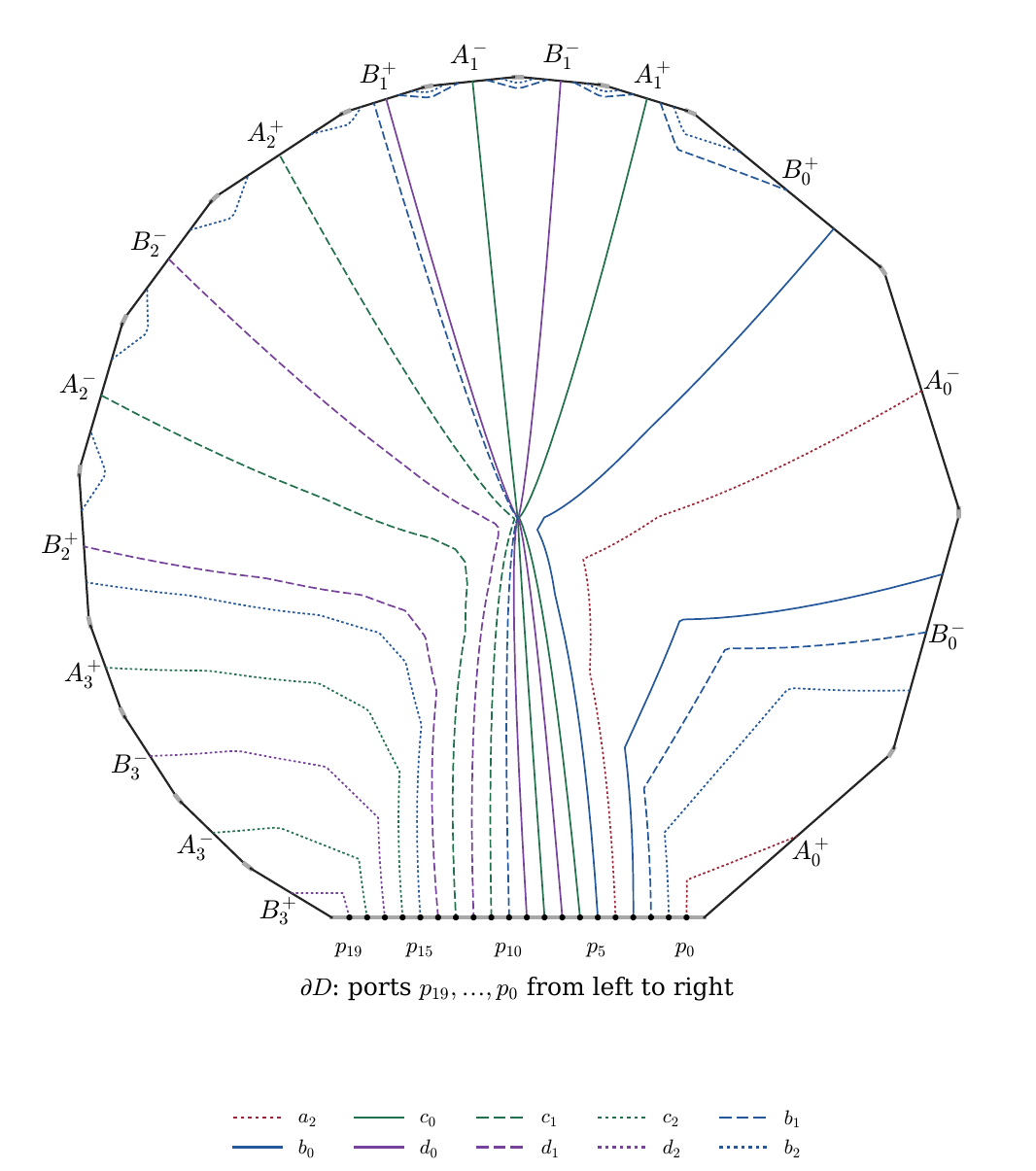}
\caption{The original non-tree edges in the paired polygon for $S$.
Thick gray intervals are parts of $\partial S$; paired labeled sides
are read in reverse order. The cap is shown separately in
Figure~\ref{fig:surfacecap}.}
\label{fig:surfacebase}
\end{figure}

Cut $S$ into a disk with thirty-two sides, alternating eight pairs of
identified sides with sixteen boundary intervals. The identified
sides have the cyclic labels
\[
 A_0^+,B_0^-,A_0^-,B_0^+,\ldots,
 A_3^+,B_3^-,A_3^-,B_3^+.
\]
Oppositely signed sides with the same letter are identified in
reverse boundary order. The sixteen intervening intervals form
$\partial S$; all twenty ports are placed on one of them.
The original graph has the routing shown in
Figure~\ref{fig:surfacebase}. Figure~\ref{fig:surfacetransport}
shows the ten non-tree edges after transport by $T$.
The retained tree is reinserted in $D$ in both cases.

\begin{theorem}\label{thm:surfacecarrier}
For $m=1$, the paired-polygon routing and the retained tree give a
simultaneous surface diagram of the transported lifted graph,
relative to $D$. Two copies of this diagram, joined through the three
relative 1-handles with the prescribed product charts, carry all
twelve relative product attaching curves $Q_{k,j}$ on an oriented
surface $H$ of genus ten. Their product framings are the surface
framings of $H$.
\end{theorem}
\begin{proof}
The routing is recorded by ordered side endpoints and their chord
pairings in \path{surface_graph_routing.json}. There are thirty-two
chords before transport and 126 after transport. The integer check
\path{check_surface_routing.py} verifies that every endpoint is
used once, that identified sides have reversed endpoint orders, and
that no two chord pairs alternate on the polygon boundary. Straight
chords in any strictly convex realization therefore have disjoint
interiors. It also reads each edge word from the side crossings,
reduces it in the free group, and checks exact agreement with the
edge table or its image under $T$, respectively. In particular this
comparison is not made only after imposing $R=1$.
The relation $T(R)=R$ holds freely as well.
The cap check \path{check_cap_tree.py} uses rational coordinates
to restore the two tree edges and the twenty port segments without
unintended intersections. It checks the eight incident directions
at every vertex against the original lifted rotation order.

These checks identify the base routing with the original marked
ribbon graph. The transported routing has the same endpoints and
relative arc homotopy classes as $T$ of that graph, using the same
specified boundary port paths. These are properly embedded arcs,
with individually fixed, distinct labeled endpoints. Embedded arc
systems with the same relative homotopy classes and fixed distinct
endpoints are isotopic relative to the boundary: remove innermost
bigons to make corresponding arcs coincide, and proceed in the
surface cut along the arcs already matched. Applying this to the ten
arcs, and reinserting the fixed tree, gives the asserted surface
isotopy. The normal ribbons are transported by this isotopy.

Take the two pages $\Sigma_4\times\{0,\pi\}$, remove the three
vertex disks in each, and join the corresponding boundary circles
through $E_0,E_1,E_2$. In a relative 1-handle, use a round
3-ball of radius $\rho$ and write
\[
 E_j=D^1_v\times D^3_{x,y,u},\qquad
 H\cap E_j=\{u=0,\ x^2+y^2=\rho^2\},\quad \rho>0.
\]
We specify the lane matching and the transverse fields at both feet.
Write $r=-j\bmod3$. The $v=-1$ foot of $E_j$ is on the
$\psi=0$ page at $V_r$, and its $v=1$ foot is on the
$\psi=\pi$ page at $V_{r-1}$. A superscript $+$ denotes an outgoing
halfedge of the oriented lifted edge, and $-$ denotes an incoming
halfedge. In the common oriented $(x,y)$ coordinates, the cyclic
order at $V_r$ is
\[
 a_r^+,\ b_r^-,\ a_{r-1}^-,\ b_r^+,\
 c_r^+,\ d_r^-,\ c_r^-,\ d_r^+.
\]
At the other foot, replace $r$ by $r-1$. These are the vertex
rotations of the retained ribbon graph. The map $T$ fixes the foot
patches, so transporting the graph does not change this order.
Choose the transverse coordinate $t$ on each oriented base band so
that $(\partial_s,\partial_t)$ gives the surface orientation.
The signed angular field at an outgoing end is then positive and
at an incoming end negative. A constant choice of sign for a product
normal vector does not change its framing class.
The following table gives the matching; all indices are modulo three,
and $\epsilon_\ell$ is the sign of $\partial_t$ relative to the
positive angular field at either foot. The labels stand for their
transported arcs under $T$.
\begin{center}
\small
\begin{tabular}{ccccc}
\toprule
Lane $\ell$ & $\psi=0$ foot & $\psi=\pi$ foot & Product curve & $\epsilon_\ell$\\\midrule
0 & $a_r^+$ & $a_{r-1}^+$ & $Q_{a,j}$ & $+1$\\
1 & $b_r^-$ & $b_{r-1}^-$ & $Q_{b,j}$ & $-1$\\
2 & $a_{r-1}^-$ & $a_{r-2}^-$ & $Q_{a,j+1}$ & $-1$\\
3 & $b_r^+$ & $b_{r-1}^+$ & $Q_{b,j}$ & $+1$\\
4 & $c_r^+$ & $c_{r-1}^+$ & $Q_{c,j}$ & $+1$\\
5 & $d_r^-$ & $d_{r-1}^-$ & $Q_{d,j}$ & $-1$\\
6 & $c_r^-$ & $c_{r-1}^-$ & $Q_{c,j}$ & $-1$\\
7 & $d_r^+$ & $d_{r-1}^+$ & $Q_{d,j}$ & $+1$\\\bottomrule
\end{tabular}
\end{center}
The two occurrences of $Q_{k,j}$ for $k=b,c,d$ are its distinct
passages through $E_j$ in opposite directions. For $k=a$, the two
passages lie in $E_j$ and $E_{j-1}$, in agreement with the attaching
word above.

Choose one oriented parametrization of the base 0-handle so that
its eight halfedges meet the boundary at
$\theta_\ell=2\pi\ell/8$ in the listed order, and lift this same
parametrization to all three vertex disks on both pages. Carry the
base-band collars and their transverse coordinates by the same
parametrization. This common choice respects the existing product
transport, which fixes the base coordinates; it does not introduce
an independent change of coordinates at either foot. The map $T$
fixes these charts. Near a lane endpoint, use the collar coordinates
\[
 (q,t,u)\longmapsto
 \bigl(R(q,u)\cos(\theta_\ell+\epsilon_\ell t/\rho),\
       R(q,u)\sin(\theta_\ell+\epsilon_\ell t/\rho),u\bigr),
 \qquad R(q,u)=\sqrt{\rho^2-u^2}+q.
\]
Here $q=0$ is the foot sphere and $u=0$ is the page; the
coordinates are restricted to a sufficiently small collar. At the
endpoint $q=t=u=0$, the derivatives send
\[
 \partial_t\longmapsto
 \frac{\epsilon_\ell}{\rho}\partial_\theta,
 \qquad \partial_u\longmapsto\partial_u.
\]
The same formulas hold at the two feet, since the fiber transport
changes the sheet index but fixes the base-band transverse
coordinate. Extend the eight lanes through the handle by
\[
 v\longmapsto(v,\rho\cos\theta_\ell,\rho\sin\theta_\ell,0),
 \qquad -1\leq v\leq1,
\]
and extend their normal fields by
$(\epsilon_\ell\partial_\theta/\rho,\partial_u)$.
Their angles and signs are constant along $v$, so the lanes are
disjoint and the extension has no annular winding. The two page
orientations used to orient $H$ do not reverse these fixed
$(x,y)$ foot coordinates.
The curve $Q_{k,j}$ uses $T(k_{-j})$ on the first page and
$T(k_{-(j+1)})$ in reverse on the second, with its other sides in
$E_{j-\delta(k)}$ and $E_j$.
Every lifted edge occurs once on each page, including the two tree
edges. Thus all twenty-four page arcs are represented together.
Orient the two pages oppositely. The annuli join them to an oriented
connected surface with
\[
 \chi(H)=2\chi(\Sigma_4)-6=-18,\qquad g(H)=10.
\]
On the pages, the product vector $\partial_t$ is tangent to the
surface and transverse to each lifted band core. The table and
collar derivatives identify this vector at both endpoints of every
lane with its specified constant angular field. Reversal of $\psi$
acts by $(x,y,u)\mapsto(x,y,-u)$ at both feet. Its derivative fixes
the displayed first normal vector and negates the second at both
ends, so it introduces no rotation in the first normal vector.
The opposite orientations of the pages agree with the opposite
traversals of the two band sides. After rounding each product corner,
the matching vectors give a continuous nonzero transverse field on
$H$ around every $Q_{k,j}$. This is the first vector of the transported
product framing. Its push-off stays on $H$, and hence the product
framing has coefficient zero relative to the surface framing.
\end{proof}

\begin{figure}[!htbp]
\centering
\includegraphics[width=.92\textwidth,height=.60\textheight,keepaspectratio]{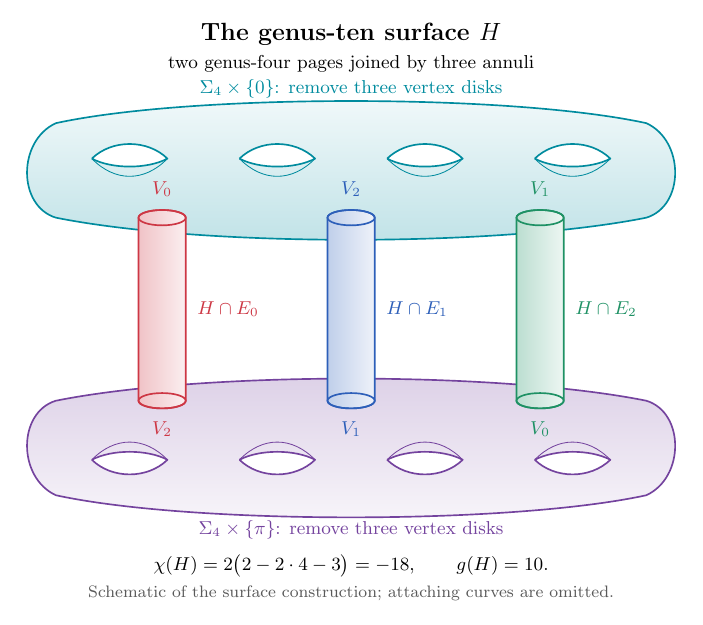}
\caption{The genus-ten surface $H$. Remove three vertex disks from each
of the pages $\Sigma_4\times\{0\}$ and $\Sigma_4\times\{\pi\}$ and
join their boundary circles by the annuli $H\cap E_j$. The annulus
through $E_j$ joins $V_{-j}$ on the first page to $V_{-(j+1)}$ on the
second, with indices modulo three. The pages are oriented oppositely.
Thus $\chi(H)=2(2-8-3)=-18$ and $g(H)=10$. This is a schematic of
the surface construction; the twelve attaching curves are not shown.}
\label{fig:genusten}
\end{figure}

Figure~\ref{fig:genusten} shows the construction of $H$ from its two
genus-four pages and the three annuli through the relative 1-handles.

\begin{figure}[!htbp]
\centering
\includegraphics[width=.68\textwidth,height=.58\textheight,keepaspectratio]{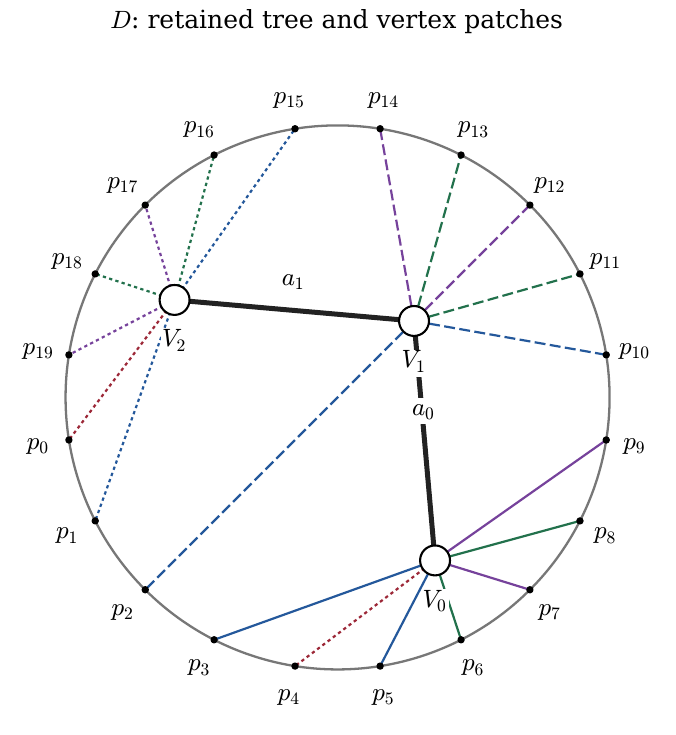}
\caption{The three vertex patches and retained tree in $D$.
Equally numbered endpoints are identified with those in the polygon
along the oppositely oriented boundaries.}
\label{fig:surfacecap}
\end{figure}

\begin{figure}[!htbp]
\centering
\includegraphics[width=.98\textwidth,height=.69\textheight,keepaspectratio]{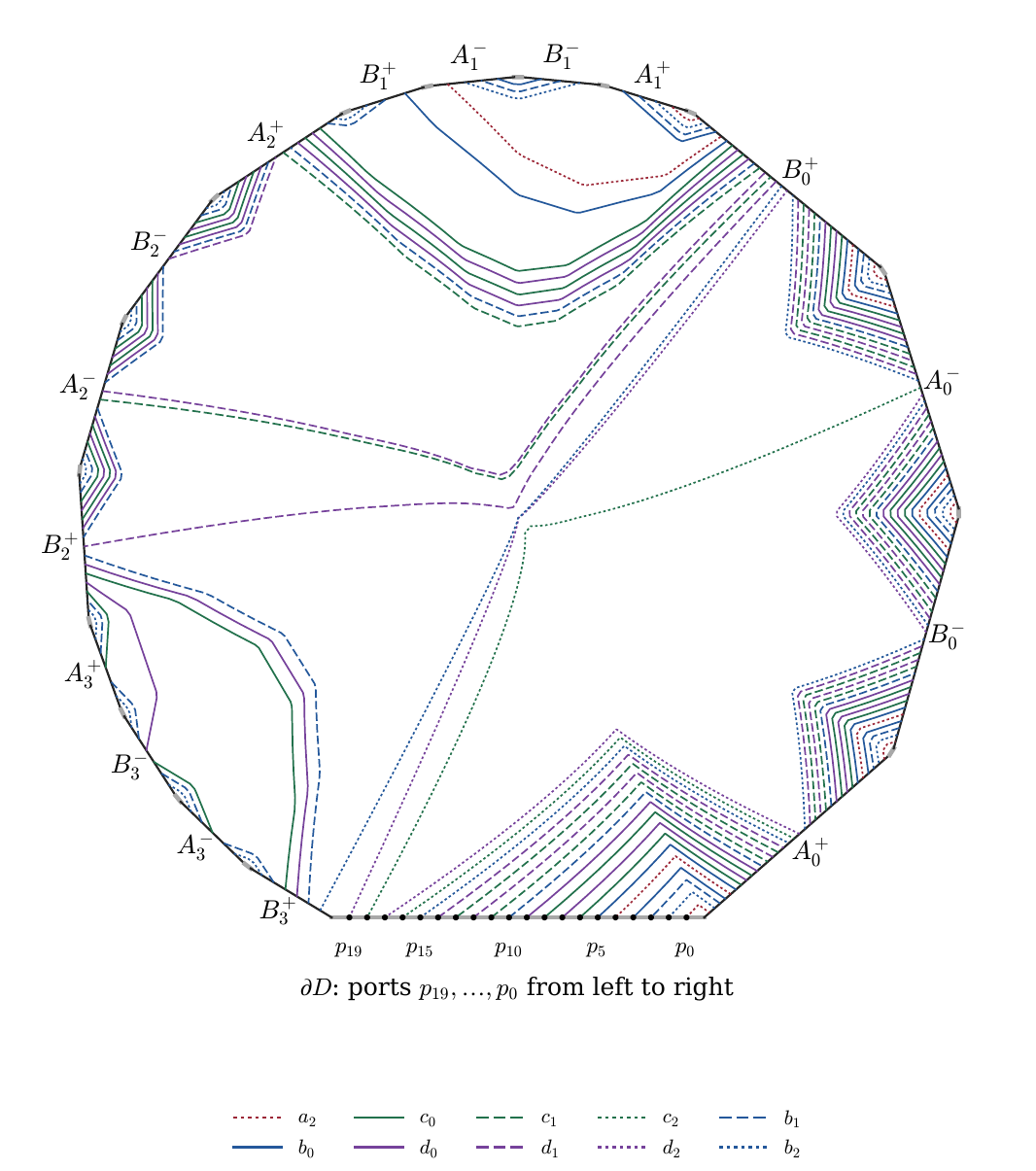}
\caption{The ten non-tree edges of the transported graph $T(\Gamma)$
for the symmetric $m=1$ sum. Each legend entry
$k_r$ denotes $T(k_r)$. Pair equally labeled sides with reversed
endpoint order, then reinsert the cap from
Figure~\ref{fig:surfacecap}. Two copies of the polygonal diagram, joined through $E_j$ along arcs
of constant angle, give the twelve relative product curves on $H$.
The zero surface framings refer to $H$. They do not assign
integer labels in a dotted-circle presentation.}
\label{fig:surfacetransport}
\end{figure}

Pullback by $(T,\psi\mapsto-\psi)$ gives the closed attaching maps;
the two cap curves retain their base-disk normal framings. Theorem
\ref{thm:surfacecarrier} gives the simultaneous product graph and
surface framings. Full handlebody diagrams of the closed sums,
with their complete attaching links and integer framing coefficients,
will be studied in the follow-up paper~\cite{AAyII}. No further
handle cancellations are made here.

\section{A complement homology check}\label{app:complementcheck}
The torsion in $H_2$ also follows from the complement exact sequence.
Put $X=\Sigma_g\times S^2$ and $F=F_{g,p}$. Excision and the Thom
isomorphism give $H_k(X,M_{g,p};\Z)\cong H_{k-2}(F;\Z)$.
The map $H_4(X)\to H_2(F)$ is an isomorphism. The relevant
part of the long exact sequence is
\[
 H_3(X)\longrightarrow H_1(F)\longrightarrow H_2(M_{g,p})
 \longrightarrow H_2(X)\longrightarrow H_0(F).
\]
The first map is Poincar\'e dual to
$q_p^*:H^1(\Sigma_g;\Z)\to H^1(F;\Z)$. The matrix $Q_{p,g}$
in Section~3 represents $q_{p*}$ and has nonzero Smith entries
$p,1,\ldots,1$. Transposition and the unimodular Poincar\'e duality
identifications preserve these entries. Thus the first map has
cokernel $\Z^{2G-2g}\oplus\Z/p$.
The intersection map
$H_2(X)\to H_0(F)$ has image $p\Z$ and kernel $\Z$, since
$[F]=p[\mathcal S]$. Hence
\[
0\longrightarrow\Z^{2G-2g}\oplus\Z/p
\longrightarrow H_2(M_{g,p};\Z)
\longrightarrow\Z\longrightarrow0.
\]
Since the quotient is free, this sequence splits as abelian groups
and recovers the stated $H_2$. This establishes the abstract group;
Corollary~\ref{cor:poletorsion} identifies the difference of the two
pole sections as a generator of its torsion subgroup.

\section{Certification of the pole projections}\label{app:crossings}
We prove Proposition~\ref{prop:certifiedcrossings} for the analytic
arcs outside the paired feet, not for the closed sums. The calculation
uses 128-bit Arb arithmetic~\cite{Arb} through \texttt{python-flint}.
The source package contains \path{certify_projection_crossings.py},
camera data in \path{figures/usual_projection_data.json}, seeds in
\path{verified_usual_projections.json}, and certificates in
\path{certified_projection_crossings.json}.
The validator \path{projection_parameters.py} checks the drawing
parameters; \path{test_projection_parameters.py} tests parameter drift.
We use $\epsilon=0.032$, $\zeta=0.012$, and $\eta=0.005$, obtained
by the framed collar deformation of Section~2. Only the center arcs
are certified, not the enlarged normal markers. Seeds choose candidate
neighborhoods; full-domain interval subdivision proves completeness.
Execution instructions and software versions for all checks are in
\path{COMPUTATIONAL_README.txt}.

The foot angles are $\theta_m=2\pi(m+0.13)/(4g)$, and the
fiber angles are $\phi=\sigma_P(2\pi j/p+0.23)$. A closing interval
has base angle between $\theta_m+\zeta$ and $\theta_{m+1}-\zeta$.
Each arc has the form $q(t)=A\cos t+B\sin t$, with $A,B$ orthonormal.
For a corner arc, take
$A=(\cos\theta,\sin\theta,0,0)$ and
$B=(0,0,\cos\phi,\sin\phi)$, with
$\epsilon\leq t\leq\pi/2-\epsilon$. For a closing interval, recenter
its base angle at $\theta_c=(\theta_m+\theta_{m+1})/2$ and use
$A=(\cos\theta_c,\sin\theta_c,0,0)$,
$B=(-\sin\theta_c,\cos\theta_c,0,0)$, with
$|t|\leq\pi/(4g)-\zeta$.
With $h=\tan(t/2)$ and $M=\mathbf E\mathbf O$, write
\[
 X_i(h)=\frac{n_i(h)}{d(h)},\qquad
 n_i(h)=(A M)_i(1-h^2)+2(B M)_ih,\quad
 d(h)=1+h^2-(A\cdot\ell)(1-h^2)-2(B\cdot\ell)h.
\]
Thus a crossing of arcs $a,b$ is a zero of the two polynomials
\[
 F_i(h,k)=n_i^a(h)d^b(k)-n_i^b(k)d^a(h),\qquad i=1,2.
\]
The script encloses each stored camera entry in a ball of radius
$10^{-12}$. It checks that the exact normalized pole and the
Gram-Schmidt orthonormalizations of the stored matrices lie in these
balls, projecting the columns of $\mathbf E$ to $\ell^\perp$ first. All subsequent tests hold uniformly over the resulting
coefficient intervals. The denominators are positive on the entire
arc domains. For each arc the coefficient matrix of
$(n_1,n_2,d)$, in the basis $(1,h,h^2)$, has determinant bounded
away from zero. Its planar map is therefore a projective image of
a nonsingular conic and is injective on the parameter interval.

For every pair of distinct arcs, subdivide the full parameter
rectangle. Discard a box if either interval evaluation of $F_i$
excludes zero. All other terminal boxes must lie in one of the
certified root neighborhoods recorded in
\path{certified_projection_crossings.json}. The entries are grouped by
chart name; each crossing records its two arc indices, root box,
root enclosure, and depth gap. For a candidate box $B$,
with center $x_0$, choose a nonsingular approximate inverse $C$ to
the Jacobian and form
\[
 K=x_0-CF(x_0)+(I-C\,DF(B))(B-x_0).
\]
The interval tests give $K\subset\operatorname{int}B$ and
$\|I-C\,DF(B)\|_\infty<1/2$. The contraction mapping theorem
then gives exactly one zero in $B$, with enclosure $K$; the same
bound makes its Jacobian nonsingular. Each root box lies strictly
inside both arc domains. Distinct root neighborhoods for one pair
are disjoint, and the projected enclosures of all crossing points
are pairwise disjoint. Finally the interval for $X_3^a-X_3^b$ at
each root excludes zero and gives the recorded overpass.

The exhaustive subdivision terminates for all six charts. The total
numbers of tested boxes are, in the chart order
$M_{2,3}^N,M_{2,3}^S,M_{3,2}^N,M_{3,2}^S,M_{2,4}^N,M_{2,4}^S$,
$67880$, $48700$, $123596$, $97272$, $94428$, and $159638$.
Lower bounds for the absolute depth gaps are respectively
$0.910$, $0.180$, $1.226$, $0.120$, $0.402$, and $0.180$.
These checks establish the stated crossing counts and overpasses
for the complement projections.

\FloatBarrier
\Needspace{7\baselineskip}

\par\vspace{-4pt}
\noindent{\footnotesize\textsc{Anar Akhmedov}\\
School of Mathematics, University of Minnesota\\
Visiting Scholar, Department of Mathematics, Harvard University\\
Email: \href{mailto:akhmedov@umn.edu}{\texttt{akhmedov@umn.edu}}\par}

\vspace{-8pt}
\begin{thebibliography}{99}\scriptsize\setlength{\itemsep}{0pt}\setlength{\parskip}{0pt}
\bibitem{AkBook} S. Akbulut, \emph{4-Manifolds}, Oxford Graduate
Texts in Mathematics, vol. 25, Oxford University Press, Oxford, 2016.
\bibitem{AkCCM} S. Akbulut, \emph{The Catanese-Ciliberto-Mendes Lopes
surface}, J. G\"okova Geom. Topol. GGT 5 (2011), 86-102.
\href{https://arxiv.org/abs/1101.3036}{arXiv:1101.3036}.
\bibitem{AkAP3} S. Akbulut, \emph{The Akhmedov-Park exotic
$\mathbb{CP}^2\#3\overline{\mathbb{CP}}^{\,2}$}, Adv. Math.
274 (2015), 928-947.
\href{https://arxiv.org/abs/1106.3924}{arXiv:1106.3924}.
\bibitem{AkAP2} S. Akbulut, \emph{The Akhmedov-Park exotic
$\mathbb{CP}^2\#2\overline{\mathbb{CP}}^{\,2}$}, Adv. Math.
274 (2015), 948-967.
\href{https://arxiv.org/abs/1107.3073}{arXiv:1107.3073}.
\bibitem{A06} A. Akhmedov, \emph{Construction of symplectic cohomology
$S^2\times S^2$}, Proceedings of the 14th G\"okova Geometry-Topology
Conference 2007, International Press, 2008, 36-48.
Preprint: \href{https://arxiv.org/abs/math/0611128}{arXiv:math/0611128v4}.
\bibitem{AB} A. Akhmedov, \emph{Braided multisections and symplectic
four-manifolds with the rational cohomology of $S^2\times S^2$}, 2026,
\href{https://arxiv.org/abs/2609.11727}{arXiv:2609.11727v1}.
\bibitem{Stein} A. Akhmedov, \emph{Braided Multisections and Aspherical
Stein Fillings}, preprint, 2026.
\bibitem{AY} A. Akhmedov and S.-K. Yeung, \emph{Braided Multisection
Sums, Cyclic Covers, and Fake Quadrics}, manuscript, 2026.
\bibitem{AAyII} A. Akhmedov, S. Akbulut, and E. Z. Yildiz,
\emph{Handlebodies for Akhmedov's Fake Quadrics}, manuscript.
\bibitem{Arb} F. Johansson, \emph{Arb: efficient arbitrary-precision
midpoint-radius interval arithmetic}, IEEE Trans. Comput. 66 (2017),
no. 8, 1281-1292. \href{https://doi.org/10.1109/TC.2017.2690633}{doi:10.1109/TC.2017.2690633}.
\bibitem{Gompf} R. E. Gompf, \emph{A new construction of symplectic
manifolds}, Ann. of Math. (2) 142 (1995), no. 3, 527-595.
\bibitem{Smith} I. Smith, \emph{Symplectic submanifolds from surface
fibrations}, Pacific J. Math. 198 (2001), no. 1, 197-205.
\bibitem{HP} C. S. Hays and B. D. Park, \emph{Non-isotopic symplectic
surfaces in product $4$-manifolds}, Trans. Amer. Math. Soc. 360 (2008),
no. 11, 5771-5788.
\end{thebibliography}
\end{document}